\documentclass[11pt]{article} 

\usepackage[utf8]{inputenc} 

\usepackage{geometry} 
\usepackage{amssymb}
\usepackage{amsthm}
\usepackage{amsmath}
\usepackage{amsfonts}

\usepackage{graphicx} 
\usepackage{subcaption}
\RequirePackage[numbers]{natbib}

\usepackage{amsmath,amssymb,amsthm,dsfont,xcolor}
\usepackage{booktabs} 
\usepackage{array} 
\usepackage{paralist} 
\usepackage{verbatim} 
\usepackage{subfig} 

\usepackage{fancyhdr} 
\usepackage{sectsty}
\allsectionsfont{\sffamily\mdseries\upshape} 

\usepackage[titles,subfigure]{tocloft} 

\newtheorem{theorem}{Theorem}
\newtheorem{lemma}{Lemma}
\newtheorem{definition}{Definition}

\newtheorem{corollary}{Corollary}

\title{Inference for concave distribution functions under measurement error}

\author{
	Mohammed Es-Salih Benjrada\\
	Department of Economics, University of Bergamo, Italy\\
	\texttt{mohammedessalih.benjrada@unibg.it}
	\and
	Cécile Durot\\
	MODAL'X, Université Paris Nanterre, France\\
	\texttt{cecile.durot@parisnanterre.fr}
	\and
	Tommaso Lando\\
	Department of Economics, University of Bergamo, Italy\\
	\texttt{tommaso.lando@unibg.it}
}

\begin{document}
\maketitle

\begin{abstract}
	We study nonparametric inference for a concave distribution function under the measurement error model, where the non-negative variable of interest is perturbed by additive independent noise. We propose a shape-constrained estimator defined as the least concave majorant on the non-negative real half-line of the deconvolution estimator of the cumulative distribution function, and we establish its uniform consistency as well as its square-root convergence in distribution. To assess the concavity assumption, we introduce a nonparametric test of the null hypothesis that the distribution function is concave on the non-negative real half-line against the alternative that it is not. The test is calibrated via the bootstrap. We show that the test statistic and its bootstrap analogue have the same limiting distribution under the null, while the rejection probability tends to one under the alternative. The proofs mainly rely on a bootstrap Donsker-type result for the deconvolution estimator of the cumulative distribution function, combined with the functional delta method. Simulation studies illustrate the finite-sample performance of both the estimator and the test.
\end{abstract}



\section{Introduction}\label{sec: intro}
Statistical inference under shape restrictions exploits nonparametric assumptions on the shape of the distribution of interest to improve estimates. Additionally, it deals with nonparametric tests about such shape constraints. The books by  \cite{robertson1988,groeneboom2014,silvapulle2011constrained} provide an overview of the issues addressed in the context of statistical inference under shape constraints. The results depend on the methodology and on the type of assumption, see for instance \cite{grenander,marshallmle,durot2008testing,rojo,landoodds,lando2024nonparametric,lando2025new}. One of the most commonly used constraints, starting from the seminal work of \cite{grenander}, is the concavity of the cumulative distribution function \citep{carolan2002,beare2021}, corresponding to a decreasing density \citep[Chapter 2.2]{groeneboom2014}.

Differently from the usual approach, in this paper we deal with this problem under the classic measurement error model (also known as deconvolution model), namely, we assume that the random variable of interest $X$ is contaminated by some additive random and independent noise $\varepsilon$, with a known distribution. In other words, instead of $X$, which we assume to be non-negative and absolutely continuous, we observe $%
Y=X+\varepsilon $. The usual way of dealing with this model is using the deconvolution estimator \citep{stefanski1990deconvolving,Meister2009}. In this paper, we propose to improve this estimator using the concavity asumption. More importantly, we will focus on testing the concavity of the cumulative distribution function of interest under the measurement error model. If we are able to confirm the concavity assumption, based on observed data, then we can use this information to improve the estimation by exploiting the shape constraint. Otherwise, if the test detects deviations from concavity, then the unconstrained estimator will be preferable. Another reason why a concavity test is particularly relevant, in the measurement error model, is that the distribution of the additive error has a crucial impact on the shape of the cumulative distribution function of {$Y$}. On the one hand, $X$ may have a concave distribution function but $Y$ may not (as a trivial example, among many others, if $X$ and $\varepsilon$ are uniform, the resulting $Y$ has a triangular distribution). On the other hand, in some cases when the distribution function of $X$ is non-concave, the convolution with $\varepsilon$ may yield a concave distribution function on $\mathbb{R}^{+}$.

This paper combines theories on deconvolution estimation (i.e., estimation under measurement error) and on estimation and testing under shape constraints.  As a shape-constrained estimator, we consider the least concave majorant of the classic kernel deconvolution estimator of the distribution function $F$ of $X$. We show its uniform consistency under general assumptions, and we establish a Donsker-type result for it, i.e. we show that the difference between the estimator and the true concave distribution function converges weakly at the rate $\sqrt n$ to a  centered Borel random variable in $\ell_\infty(\mathbb R)$. This latter result exploits the results of \cite{sohl2012uniform} combined with the functional delta method.

Then, in order to test concavity of $F$, we consider the test statistic defined as the supremum distance between the kernel deconvolution estimator and its least concave majorant. 
A similar approach has been developed in several settings to test concavity of a given function (or monotonicity of its derivative when it exists). Depending on the setting, the function of interest could be, for instance, a cumulative hazard rate, a distribution function, a copula function. The general idea is to define the test statistic as the $L_p$-distance, with some fixed $p\in[1,\infty]$, between an unconstrained estimator of the target function and its least concave majorant.
Two different strategies have been proposed to calibrate such a test. The first one relies on a least favorable hypothesis, as discussed in \cite{durot2003kolmogorov,durot2008testing, kulikov2004testing, carolan2005nonparametric,
	delgado2012distribution,
	beare2015nonparametric}. This approach consists in showing that a linear target function is least favorable in the sense that the test statistic is (at least asymptotically)
stochastically larger under linearity than it would be under a concave but nonlinear function. The limiting distribution is then derived under this least favorable case, and the  corresponding $(1-\gamma)$-quantile  is used as the critical value, ensuring that the test achieves the desired asymptotic significance level  $\gamma$. The second strategy consists in deriving the limit distribution of the test statistic under every  arbitrary concave functions and to approximate the limit using either plug-in or bootstrap. This gives more powerful tests than those based on a least favorable hypothesis. However, a major difficulty with this approach stems from the fact that the rate of convergence of the test statistic differs depending on whether or not the target function is strictly concave. For instance in case of a distribution function on $[0,1]$, it is of order
$n^{1/2}$ under linearity and $n^{2/3}$ under strict concavity with appropriate smoothness conditions, see \cite{kulikov2004testing}. 
The latter paper points out that, due to different rates of convergence, scaling the test statistic by $n^{2/3}$ necessarily leads to a test that does not have the prescribed asymptotic level. On the other hand, scaling the test statistic by $n^{1/2}$ leads to  a degenerate limit distribution under strict concavity with no affine segments. This makes it difficult to prove the consistency of a bootstrap calibration since such proofs rely on continuity of the limit distribution, see e.g. \cite[Chapter 23]{van2000asymptotic}. Moreover, when non degenerate, the limiting distribution is too complex to be feasibly approximated using plug-in methods. In three different settings, \cite{seo2018tests, beare2019improved,fang2019refinements} develop such a methodology, and prove that their test statistic and its bootstrap version share the same limit distribution under the null hypothesis. Hence, the asymptotic probability of rejection is the prescribed level  under each concave function for which the limit distribution of the corresponding test statistic is continuous. It is unclear what is the asymptotic probability under strict concavity with no affine parts since the limit distribution degenerates at zero in that case. \cite{seo2018tests, beare2019improved} have neither simulations nor theoretical results for that case. On the other hand, \cite{fang2019refinements} slightly modifies the test statistic to cover some cases of strict concavity albeit some other cases are not covered, and simulations show good performance even under some cases of strict concavity. 

{In the model with measurement error,  we do not have access to a least favorable hypothesis (see Section \ref{sec: lfh}) so we scale the test statistic by $n^{1/2}$ and consider a bootstrap calibration. In the same spirit as in \cite{seo2018tests, beare2019improved,fang2019refinements},
	computing the limit distribution of the rescaled test statistic under arbitrary concave function relies on computing the limit distribution of the unconstrained estimator (provided in our setting by  \cite{sohl2012uniform}), and combine it with Hadamard directional differentiability of the least concave majorant operator and with the delta method. We build a boostrap statistic that differs in spirit from those in \cite{seo2018tests, beare2019improved,fang2019refinements}, and that has the same limit distribution as the original one under any concave distribution function. The bootstraps used in the three papers we have just cited use a similar strategy, inspired by \cite{ dumbgen1993nondifferentiable,
		fang2019inference},  to circumvent the fact that the standard bootstrap does not work, due to the very weak form of differentiability of the least concave majorant operator. For our part, we circumvent this difficulty by \lq\lq sub-sampling\rq\rq\
	in the sense that our bootstrap sample has size 
	$m\ll n$. Moreover, whereas the above papers sample from a sort of empirical distribution and then adapt their bootstrap test statistic so that it mimics  the behaviour the original one would have under the null hypothesis, we directly sample from a distribution that satisfies the null hypothesis. We show that if the distribution function is concave but not strictly concave (meaning that it has at least an affine part), the probability of rejection converges to the prescribed level  as the sample size goes to infinity. On the other hand, the common limit distribution is again degenerate at zero under strict concavity with no affine segments. We provide some heuristics on the fact that the asymptotic probability of rejection in that case is smaller than the prescribed level, and we show simulation experiments in which the probability becomes very close to zero as the sample size increases. 
}

The performance of the constrained estimator and that of the test are illustrated through a simulation study. The simulations suggest that, for sufficiently large sample sizes and a moderate noise-to-signal ratio, (i) when $F$ is concave, the constrained estimator outperforms the unconstrained one; (ii) under concavity of $F$, the rejection rate of the test does not exceed the precribed level; and (iii) the test can detect nonconcavity of $F$ even in settings where the distribution of $Y$ is concave.

As by-products of our proofs, we provide 
\begin{enumerate}
	\item
	a bootstrap version of  the Donsker Theorem 1 in  \cite{sohl2012uniform}  from which, as a particular case, one can derive the limit behavior of the deconvolution estimator of the distribution function in a bootstrap setting, under measurement error, see Theorem \ref{theo: donsker} in Section \ref{sec: generalization} and \eqref{eq: cvboot} below;
	\item
	a general result on convergence of least concave majorants (and their slopes) of a sequence of continuously differentiable functions whose sequence of slopes uniformly converge, see Theorem \ref{theo: cvunif} below.
\end{enumerate}

The paper is organized as follows. Section \ref{sec: context} provides the model, definitions and assumptions.  The constrained estimator is studied in Section \ref{sec: concaveestim}, which also contains new results about uniform convergence of least concave majorants of continuously differentiable functions, as well as that of their slopes.
The test statistic is defined in Section \ref{sec: teststat}, where its limit behavior is studied under both the null and the alternative hypotheses.
A bootstrap calibration of the test is given in Section \ref{sec: calib}, which also contains an alternative calibration that provides a less powerfull test, compared to the bootstrap, but is easier to implement. The main part of the paper ends with Section \ref{sec: simul}, which is devoted to simulations.
The main sections are followed by appendices that are organized as follows. 
Section \ref{sec: generalization} provides a generalization of the Donsker result in  \cite{sohl2012uniform} to a bootstrap setting.
The results from Sections \ref{sec: context} and \ref{sec: concaveestim}  are proved in Section \ref{sec: proofcvunif}.
Theorem \ref{theo: limit}, that provides the limit behavior of the test statistic,  is proved in \ref{sec: prooflimit}, while the proofs of Theorem \ref{theo: bootstrap} and Theorem \ref{lem: level}, that prove the validity of bootstrap calibration, are given in Section \ref{sec: proofbootstrap}. The results from  Section \ref{sec: generalization} are proved in Section \ref{sec: proofdonsker}.

\section{Context}\label{sec: context}
\subsection{Notations}
Let $\mathbb{R}^+=[0,\infty)$.

The convolution operator is denoted by $\star$.

For all $x\in\mathbb R$, $\langle x\rangle=(1+x^2)^{1/2}$. 

For all intervals $I\subset\mathbb{R}$, $\ell_\infty(I)$ denotes the space of all bounded  functions $\theta:I\to\mathbb{R}$ equipped with the sup-norm $\|\theta\|_{I\infty}=\sup_{t\in I}|\theta(t)|$, where the indice $I$ will be omitted when $I=\mathbb{R}^+$ or when there is no possible confusion; and $C^{\infty }(I)$ denotes the set of infinitely differentiable functions on $I$. For \( p \in [1, \infty) \), \( L^p(I) \) denotes the space of real-valued functions $g$ on  \( I \) such that $ \int_I |g(x)|^p \, dx $ is finite, equipped with the \( L^p \)-norm \( \|\cdot\|_p \) such that for all $g$ in the space, $\|g\|_p^p$ is equal to the latter finite integral. 

Weak convergence is denoted by \( \rightsquigarrow \) and convergence in distribution of random variables is denoted with $ \overset{L}{\longrightarrow }$. 

The probability density functions of $X$, $Y$ and $\varepsilon$ are denoted by $f$, $f_Y$ and $f_\varepsilon$, the corresponding cumulative distribution functions are denoted by $F$, $F_Y$ and $F_\varepsilon$ 	and the corresponding Fourier transforms of  are denoted by $\varphi$, $\varphi_Y$ and $\varphi_\varepsilon$, respectively. Hence in particular, $\varphi_\varepsilon=\mathcal F f_\varepsilon$, and similarly for the other two distributions.

We write $m\ll n$ to mean that $m$ depends on $n$ and is such that $mn^{-1}$ converges to zero as $n\to\infty$.

\subsection{Model, objective and basic estimators}\label{sec: model}

Under measurement error, we observe a random sample $Y_1,\dots,Y_n$ from a random variable $Y=X+\varepsilon$, where $X$ is a non-negative random variable of interest and $\varepsilon$ is a random noise, which we assume to be independent from $X$. In particular, the sample can be written in the form $Y_{i}=X_{i}+\varepsilon _{i}$, for $i=1,\dots,n,$ where $X_i$ and $\varepsilon_i$ are iid copies of $X$ and $\varepsilon$, respectively. 

The distribution of $X$ is not identifiable if that of $\varepsilon$ is unknown. Consequently, assuming a form for the distribution of $\varepsilon$ is frequent \citep{Fan1991a,stefanski1990deconvolving,adusumilli2020inference}. We assume that the random noise $\varepsilon$ have a known probability density function, denoted with $f_\varepsilon$. Moreover, we assume that $X$ has a probability  density function as well. We denote the cumulative distribution function and the probability density function of $X$ with $F$ and   $f$, respectively. Then, $Y$ has a probability density function  that we  denote by $f_Y$, and that is the convolution between $f$ and $f_\varepsilon$: for all $t\in\mathbb{R}$,
\begin{equation*}
	f_Y(t)=f\star f_\epsilon (t)=\int f(t-x)f_\varepsilon(x)dx.
\end{equation*}

Our aim is to obtain a test for the null hypothesis 
\begin{equation}\label{eq: H0}
	\mathcal{H}_{0}:\ \lq\lq F\text{ is concave on }[0,\infty)",
\end{equation} against the alternative hypothesis \( \mathcal{H}_{1}:\ \lq\lq F\text{ is not concave on }[0,\infty)\)". 
Moreover, if the null hypothesis is not rejected, we wish to construct an estimator of $F$ that takes it into account, i.e. a concave estimator. Defining the concave estimator and the test requires to first define Fourier transforms, which we will now do. 

For an absolutely integrable function $g$, we define its corresponding Fourier transform as $\mathcal{F}g(t) = \int e^{its} g(s) \, ds$ for all $t\in\mathbb{R}$; note that if $Z$ is a random variable with density function $g$, then $\mathcal{F}g(t)$ represents the characteristic function $E(\exp(itZ))$ evaluated at $t$. 
We will estimate the density function $f$  using the so-called deconvolution
approach (see, e.g., \cite{stefanski1990deconvolving}), and  we obtain
an estimator of $F$ by integration. The deconvolution estimator is constructed by observing that, under the deconvolution model, we have $\mathcal{F}f_{Y} = \mathcal{F}f \times  \mathcal{F}f_{\varepsilon}$, so, if $\mathcal{F}f_{\varepsilon }$ has no vanishing
points on the real line (that is, $ \mathcal{F}f_{\varepsilon}(t) \neq 0$ for all $t\in\mathbb{R}$), then  $\mathcal{F}f = \mathcal{F}f_{Y}/\mathcal{F}f_{\varepsilon}$, and by applying the inverse Fourier transform $\mathcal{F}^{-1}$, we obtain
\begin{equation*}
	f( x) =\mathcal{F}^{-1}\left[ \frac{\mathcal{F}f_{Y}}{\mathcal{F}
		f_{\varepsilon }}\right] ( x)
\end{equation*}
for all $x\in\mathbb{R}$.
An empirical version of $\mathcal{F}f_{Y}$ can be obtained using the plug-in
method combined with the kernel method, i.e., $\varphi
_{n}\times \mathcal{F}K_{h_n}$, where $\varphi _{n}(t):=n^{-1}
\sum_{i=1}^{n}\exp (itY_{i})$ and $K_{h_n}(x)=h_n^{-1}K\left( x/h_n%
\right) $, with $K$ being a kernel and $h_n>0$ representing the bandwidth.
Consequently, the deconvolution estimator of $f$ is expressed as follows:
\begin{equation}\label{fn}
	f_{n}^{d}\left( x\right) :=\mathcal{F}^{-1}\left[ \mathcal{F}K_{h_n}\frac{\varphi
		_{n}}{\mathcal{F}f_{\varepsilon }}\right] \left( x\right)
\end{equation}
for all $x\in\mathbb{R}$. 
Note that $f_n^d$ is finite and well defined under the usual assumptions that $\mathcal{F}K$
is compactly supported and $\mathcal{F}f_{\varepsilon }$ has no vanishing
points on the real line. By integration, we get the basic estimator of
$F$ as follow:%
\begin{equation}\label{Fn}
	F_{n}^{d}\left( x\right) :=\int_{0 }^{x}\mathcal{F}^{-1}\left[ 
	\mathcal{F}K_{h_n}\frac{\varphi _{n}}{\mathcal{F}f_{\varepsilon }}\right] \left(
	t\right) dt
\end{equation}
for all $x\geq 0$. We say that it is a basic estimator because it does not satisfy a concavity constraint (concave estimators will be considered in the Section \ref{sec: concaveestim}). By definition, the estimator takes value zero at point zero, however, it does not necessarily tend to one as a distribution function should. Hence, in some cases, we will prefer to use a normalized version of $F_n^d$, defined by
\begin{eqnarray}\label{eq: genuineDF}
	\hat F_n^d(t)=\frac{F_n^d(t)}{\lim_{t\to\infty}F_n^d(t)},
\end{eqnarray}
where the denominator is well defined and positive thanks to Lemma \ref{lem: welldef} below. Note that the assumptions of the lemma (which implicitely assume that $\varphi_\epsilon$ has no vanishing point) are satisfied under our assumptions {\bf K}, {\bf D} and {\bf E}, but are much more general.

\begin{lemma}\label{lem: welldef}
	Let $f_n^d$ and $F_n^d$  be defined by \eqref{fn} and \eqref{Fn}. Assume that $h_n^{ 2+2\beta }n\rightarrow \infty $,  $\mathcal F K$ is supported on $[-1,1]$, there exist positive $\beta$ and $C$ such that $|1/\varphi_\varepsilon(t)|\leq C(1+|t|^2)^{\beta/2}$ for all $t\in\mathbb R$, $f$ is bounded and  $\int|\mathcal F f|$ is finite. We then have
	\begin{enumerate}
		\item
		$
		\|f_n^d-f\|_\infty=o_p(1)
		$ and 
		$
		\|F_n^d-F\|_\infty=o_p(1)
		$. 
		\item
		The denominator in \eqref{eq: genuineDF} tends to one in probability; hence it is well defined, finite and positive with probability that tends to one as $n\to\infty$. 
	\end{enumerate}
\end{lemma}

\subsection{Assumptions}\label{sec: assumptions}

A key ingredient in our proofs is Theorem 1 in \cite{sohl2012uniform}. This is a Donsker-type result that gives, as a special case, the limiting behavior of the deconvolution estimator of the distribution function in the measurement error model. Assumptions {\bf K}, {\bf D} and {\bf E} below are taken from that theorem.

Concerning the kernel $K$, we will need the following assumptions. Note that the deconvolution estimator $f_n^d$ takes real values if $K$ is symmetric and real valued, two assumptions that we assume to hold. The assumptions below are satisfied for instance by the kernel $K$ defined by \eqref{eq:kernel} below
where \( r \geq 2 \) is an even integer and \( s \geq 1 \) is an integer, see the third claim in Lemma \ref{lem: DE}  below.\\

\noindent
{\bf Assumptions K}
{\it 
	\begin{enumerate}
		\item$K\in L^{1}(
		\mathbb{R}) \cap \ell_{\infty }(\mathbb{R})$ is symmetric and band-limited with support of $\mathcal{F}K$ contained in $\left[ -1,1\right] .$
		
		\item  $\int K=1$, $
		\int x^{j}K\left( x\right) dx=0$ and $\int \vert x^{L+1}K\left(
		x\right)\vert dx<\infty,$ for $j=1,...,L$ and some integer $L>0$.
		
		\item $K$ has a continuous derivative and, for a positive constant $C$, it satisfies $
		\vert K(x) \vert +\vert K'(x)\vert \leq C ( 1+x^2) ^{-1}.$
	\end{enumerate}
}

On the other hand, we will need assumptions on the distribution of $X$. The probability density function $f$ of $X$ is assumed to satisfy the following
assumptions:\\

\noindent
{\bf Assumptions D}
{\it 
	\begin{enumerate}
		\item  $f$  is  bounded and satisfies $\int
		\left\vert x\right\vert ^{2+\delta }f\left( x\right) dx<\infty $ for some
		positive $\delta $.
		
		\item  $f\in H^{\alpha }\left( \mathbb{R}\right) $, the Sobolovev space of order $\alpha \geq 0$ defined by
		\begin{equation*}
			H^{\alpha }\left( \mathbb{R}\right)=\left\{f\in L_2(\mathbb{R})\ s.t.\ \int\left(1+x^2\right)^{\alpha}|\mathcal F f(x)|^2dx<\infty\right\}.
		\end{equation*}
	\end{enumerate}
}

Finally, we make assumptions on the distribution of the error. The decay rate of the characteristic function of the error ${\varepsilon}$ reflects the smoothness of the error. Here, we consider the case of ordinary smooth error, which means that the characteristic function of $\varepsilon$ decays with polynomial rate. The conditions are similar to the
classical decay condition introduced by \cite{Fan1991a} and are satisfied by many
ordinary smooth error distributions, such as the gamma or the Laplace. Specifically, we
assume that:\\

\noindent
{\bf Assumptions E}
{\it 
	\begin{enumerate}
		\item$\int \left\vert x\right\vert ^{2+\delta }f_{\varepsilon
		}\left( x\right) dx<\infty $ for some positive $\delta $.
		
		\item  $\varphi_{\varepsilon }$ has no vanishing
		points on the real line, and for some \( \beta > 0 \) and $C_1>0$, the following condition holds:
		\begin{equation}\label{eq: hypphiepsion}
			\left\vert \left( 1/\varphi_{\varepsilon} \right)^{\prime}(t) \right\vert \leq C_1 (1 + |t|^2)^{(\beta - 1)/2}, \quad \text{for all } t \in \mathbb{R}.
		\end{equation}
	\end{enumerate}
}

\section{Concave estimators}
\label{sec: concaveestim}

\subsection{Definition of the concave estimator and consistency}\label{sec: defestim}
A natural way of estimating $F$, exploiting the concavity assumption, is using the least concave majorant, defined below, of a basic estimator.

\begin{definition}\label{def: LCM}
	\label{def1} Given an interval $I\subseteq \mathbb{R}^{+}$, the operator  $\mathcal{M}_{I}:\ell^{\infty }\left(
	I\right) \rightarrow \ell^{\infty }\left( I\right) $ maps each function $%
	\theta \in \ell^{\infty }\left( I\right) $ to $\mathcal{M}_{I}\theta$, the least concave majorant of $\theta$ on $I$, defined for $x\in I$ by
	\begin{equation*}
		\mathcal{M}_{I}\theta \left( x\right) :=\inf \left\{ g\left( x\right) :g\in \ell^{\infty}(I) 
		\text{ is concave and }
		\theta \leq g\text{ on }I\right\}.
	\end{equation*}
	For convenience, we write $\mathcal{M}$ instead of $\mathcal{M}_{ \mathbb{R}^{+}}$ and we call this operator \lq\lq least concave majorant" without specifying the interval on which it is taken. 
\end{definition}

The following theorem, which is proved in Section \ref{sec: proofcvunif}, provides a general result on convergence of least concave majorants (and their slopes) of a sequence of continuously differentiable functions whose sequence of slopes uniformly converge.

\begin{theorem}\label{theo: cvunif}
	Let $f$ and $f_n$, for positive integers $n$, be continuous functions from $[0,\infty)$ to $\mathbb R$ and let $F_n(t)=\int_0^tf_n(u)du$ and $F(t)=\int_0^tf(u)du$ for 
	$t\geq 0$. Assume $F_n,F\in\ell_\infty(\mathbb R^+)$ and denote by $f_{n0}$ and $f_0$ the respective slopes of  $\mathcal MF_n$ and $\mathcal MF$.  Then one has:
	\begin{enumerate}
		\item
		$f_{n0}$ and $f_0$ are continuous and non-negative. 
		\item If moreover $f$ is non-negative and $
		\|f_n-f\|_\infty=o(1)$ as $n\to\infty$ then 
		$
		\|\mathcal M F_n-\mathcal M F\|_\infty=o(1)
		$ and $
		\|f_{n0}- f_0\|_\infty=o(1)
		$.
		\item In a stochastic setting where $f_n$ and $f$ are possibly random, if $f$ is non-negative and $\|f_n-f\|_\infty=o_p(1)$ as $n\to\infty$,  then $
		\|\mathcal M F_n-\mathcal M F\|_\infty=o_p(1)
		$ and $
		\| f_{n0}-f_0\|_\infty=o_p(1) 
		$. 
	\end{enumerate}
\end{theorem}
As a concave estimator of $F$, we will consider the least concave majorant $\mathcal M F_n^d$ with $F_n^d$ denoting the basic estimator taken from \eqref{Fn}. In some cases, we will prefer to use the least concave majorant of the normalized version $\hat F_n^d$  taken from \eqref{eq: genuineDF}. In particular, $\hat{F}_n^d$ seems to be more convienent for estimating $F$, while we will use $F_n^d$ for the test.
Note that both estimators $\mathcal M F_n^d$ and $\mathcal M \hat F_n^d$ are continuous by concavity and that according to \eqref{Fn} and Definition \ref{def: LCM} they take the constant value 0 on the negative half-real line. Moreover, it follows from usual properties of least concave majorants (see for instance Section 2 in \cite{durot2003distance}) that the two estimators are connected through
\begin{eqnarray}\label{eq: genuineDF2}
	\mathcal M \hat F_n^d(t)=\frac{\mathcal M F_n^d(t)}{\lim_{t\to\infty}F_n^d(t)}.
\end{eqnarray}

As a corollary (also proved in Section  \ref{sec: proofcvunif}) of Theorem \ref{theo: cvunif}, we derive uniform convergence of our concave estimators.

\begin{corollary}\label{cor}
	Let $f_n^d$ and $F_n^d$  be defined by \eqref{fn} and \eqref{Fn} in the deconvolution model of Section \ref{sec: model}. Assume that $h_n^{ 2+2\beta }n\rightarrow \infty $,  $\mathcal F K$ is supported on $[-1,1]$, there exist positive $\beta$ and $C$ such that $|1/\varphi_\varepsilon(t)|\leq C(1+|t|^2)^{\beta/2}$ for all $t\in\mathbb R$, $f$ is bounded continuous and  $\int|\mathcal F f|$ is finite. Then 
	$$\|\mathcal M _n-\mathcal MF\|_\infty=o_p(1)
	\quad and \quad
	\|f_{n0}-f_0\|_\infty=o_p(1)
	$$ where $f_0$ is the slope of $\mathcal M F$,   $\mathcal M_n$ is either $\mathcal M F_n^d$ or $\mathcal M \hat F_n^d$, and $ f_{n0}$ is  the slope of $\mathcal M_n$.
\end{corollary}

\subsection{Limit distribution}
Next, we are interested in computing the limit distribution.
A key ingredient for that is Theorem 1 of \cite{sohl2012uniform}, that provides a uniform central limit theorem for kernel estimators  in the deconvolution model. The theorem covers translation classes of linear functionals of the density of $X$ of the form $t\mapsto \vartheta_t:=\int\zeta(x-t)f(x)dx$, where the special case $\zeta:=\mathds{1}_{(-\infty,0]}$ leads to the estimation of the distribution function $\vartheta_t=F(t)$ with an estimator that ressembles that in \eqref{Fn}. Precisely, the estimator here is
\begin{equation}\label{tildeFn}
	\tilde F_{n}^{d}\left( x\right) :=\int_{-\infty}^{x}\mathcal{F}^{-1}\left[ 
	\mathcal{F}K_{h_n}\frac{\varphi _{n}}{\mathcal{F}f_{\varepsilon }}\right] \left(
	t\right) dt
\end{equation}
which is connected to the estimator in \eqref{Fn} by $F_n^d=\tilde F_n^d-\tilde F_n^d(0)$ on $[0,\infty)$. The theorem is a Donsker type of theorem in the sense that it proves that the difference between the process $(\vartheta_t)_t$ and its estimator converges in law at the rate $\sqrt n$ to a centered Gaussian Borel random variable in $\ell_\infty(\mathbb{R})$. Our assumptions {\bf K}, {\bf D} and {\bf E} are taken from that theorem.
Combining Theorem 1 of \cite{sohl2012uniform} with  Example 1 of the same paper proves the simplified following version of the theorem, for the special case $\vartheta_t=F(t)$.

\begin{theorem}[\cite{sohl2012uniform}]\label{theo: ST}
	Assume conditions \textbf{D} and \textbf{E}, with $\beta<1/2$, $\alpha +3\gamma >2\beta +1$ for some $\gamma \in(\beta, 1/2)$ and \textbf{K} with $L = \left\lfloor \alpha + \gamma \right\rfloor$. Furthermore, let $h_{n}^{2\alpha +2\gamma }n\rightarrow 0$ and  $h_{n}^{\rho }n\rightarrow \infty $ for some $\rho
	>4\beta -4\gamma +2$. Then, we have 
	\begin{equation}\label{eq: ST}
		\sqrt{n}( \tilde F_{n}^{d}-F) \rightsquigarrow \mathbb{G}_{F}\ \ in\ \ell_\infty(\mathbb{R}),
	\end{equation}%
	as $n \to \infty$, where $\mathbb{G}_F$ is a centered Gaussian Borel random element in $\ell_\infty(\mathbb{R})$, with covariance function given by
	
	$$
	\Sigma_{s,t} = \int g_t(x) g_s(x) f_Y(x) \, dx - F(t)F(s),
	$$
	where $s,t\in\mathbb{R}$ and
	
	$$
	g_t = \mathcal{F}^{-1} \left[ \frac{1}{\mathcal{F}f_\varepsilon(-\bullet)} \right] \star \mathds{1}_{(-\infty, t]}.
	$$
\end{theorem}

The limit distribution of the estimator $\mathcal MF_n ^d$ from Section \ref{sec: defestim} can be derived from the previous theorem combined with the functional delta method by applying the functional $\mathcal {M}$. 
The standard version of the functional delta method assumes the operator is Hadamard differentiable (see, e.g., Theorem 2.1 in \cite{aadbook}), but the approach is in fact broadly applicable under Hadamard directional differentiability, see Theorem 2.1 in \cite{shapiro1991asymptotic}. Hence, a key ingredient for our purpose  is Proposition 2.1 in  \cite{Brendan2017}, according to which $\mathcal M$ 
is Hadamard directionally differentiable with a derivative in closed form. For completeness, the proposition is recalled in Theorem \ref{theo: brendan} below, where the following notations are used. For any concave function $\theta\in\ell_\infty(\mathbb{R}^+)$ and $x\geq 0$, let $I_{\theta,x}$ be the union of the singleton $\{x\}$ and all open intervals $A\subset\mathbb{R}^+$ such that $x\in A$ and $\theta$ is affine on $A$. Note that $I_{\theta,x}=\{x\}$ if $\theta$ is not affine in a neighborhood of $x$, whereas it is an open interval with non-empty interior otherwise. Moreover, we use the notations $\mathcal M$ and $\mathcal{M}_{I}$ for intervals $I$ as in Definition \ref{def: LCM}, and we denote by $C_0$ the collection of continuous real valued functions on $\mathbb R^+$ vanishing at infinity, equipped with the sup-norm.

\begin{theorem}[\cite{Brendan2017}, Proposition 2.1]
	\label{theo: brendan} 
	The operator $\mathcal M$ is Hadamar 
	directionally differentiable at any concave function $\theta\in\ell_\infty(\mathbb{R}^+)$ tangentially to $C_0$. Its directional derivative $\mathcal M_\theta':C_0\to\ell_\infty(\mathbb{R}^+)$ is uniquely determinated as follows: for any $h\in C_0$, and $x\geq 0$, we have  $\mathcal M_\theta'h(x)=\mathcal M_{I_{\theta,x}}h(x)$.
\end{theorem}

In particular, if $\theta$  is not affine in a neighborhood of $x$, then $\mathcal M_\theta'h(x)=h(x)$ for all $h\in C_0$.

We are now in position to state the main result of the section, which gives the limit distribution of the concave estimator $\mathcal MF_n^d$. The theorem shows that the estimator converges to the least concave majorant of $F$ at the rate $\sqrt n$, so, in particular, the estimator converges to $F$ if the latter is concave on its support. More precisely, the difference between the estimator and the least concave majorant of $F$ converges at the rate $\sqrt n$ to a centered  Borel variable in $\ell_\infty(\mathbb R)$.

\begin{theorem}\label{theo: limitdistrib}
	Let the assumptions of Theorem \ref{theo: ST} be satisfied.  Then as $n\to\infty$, one has
	\begin{equation}\label{eq: limitdistrib}
		\sqrt n (\mathcal MF_n^d-\mathcal M F) \rightsquigarrow  \mathcal{M}
		_{F}^{\prime }\mathbb{G}_F-\mathbb{G}_F(0)\quad in\ \ell_\infty(\mathbb R^+).
	\end{equation}
\end{theorem}

Using that the least concave majorant of an arbitrary continuous function $G\in\ell_\infty(\mathbb R^+)$ coincides with $G$ at point $0$, one can see that  the definition of $\mathcal{M}
_{F}^{\prime }$ implies that $\mathcal{M}
_{F}^{\prime }G(0)=G(0)$. In particular, the right-hand process in \eqref{eq: limitdistrib} (as well as the left-hand process) takes value 0 at point 0.

\section{The test statistic}\label{sec: teststat}
To test  the null hypothesis \eqref{eq: H0} against the alternative that  $F$ is not concave on $\mathbb{R}^+$, we consider a test statistic defined as a distance between the general deconvolution estimator $F_n^d$ (that performs well under both the null and the alternative hypotheses) and an estimator of $F$ that works well under the null hypothesis only, i.e. when $F$ is concave. For the latter, we consider the least concave majorant $\mathcal MF_n^d$ of $F_n^d$  (see Definition \ref{def: LCM}) which, according to Corollary \ref{cor}, uniformly converges to $\mathcal MF$, which is precisely $F$ under the null hypothesis, and differs from $F$ otherwise.  Hence, the test statistic typically takes small values under the null hypothesis and large values under the alternative hypothesis, so the null hypothesis will be rejected if the test statistic is too large.  The test statistic is formally defined as
\begin{eqnarray}
	T_n( F_{n}^{d}) &:=&\sqrt{n}\sup_{t\geq 0} \left\vert\mathcal MF_n^d(t) -F_{n}^{d}(t)\right\vert 
	\label{Tn}
	\\ \notag
	&=&\sqrt{n}\Vert\mathcal MF_n^d -F_{n}^{d}\Vert_\infty .
\end{eqnarray}

Note that due to \eqref{eq: genuineDF2} and \eqref{eq: lim F_n^d} below, one has $T_n( \hat F_{n}^{d}) =T_n( F_{n}^{d}) (1+o_p(1))$ so replacing $ F_{n}^{d}$ by $\hat F_{n}^{d}$ would provide a test with similar performances. Hence, for simplicity, we consider the test based on the non-normalized estimator $F_n^d$. 

In order to calibrate the test, i.e. to give a precise meaning of the words \lq\lq too large" above, we need to compute the limiting distribution of the test statistic under the null hypothesis.  This is the aim of the following theorem, which also gives its general behavior under the alternative hypothesis. Note that under the alternative hypothesis, we are able to consider very mild assumptions since we do not need to compute the limiting distribution of the test statistic.

\begin{theorem}\label{theo: limit}
	\begin{enumerate}
		\item
		Let the assumptions of Theorem \ref{theo: ST} be satisfied.  If  $F$ is concave on $\mathbb{R}^+$, then as $n\to\infty$, one has
		\begin{equation}\label{convergence}
			T_n( F_{n}^{d}) \overset{L}{\longrightarrow }\sup_{t\geq 0} \left\vert \mathcal{M}
			_{F}^{\prime }\mathbb{G}_F(t)-\mathbb{G}_F(t)\right\vert .
		\end{equation}
		\item Assume that $h_n^{ 2+2\beta }n\rightarrow \infty $,  $\mathcal F K$ is supported on $[-1,1]$, there exist positive $\beta$ and $C$ such that $|1/\varphi_\varepsilon(t)|\leq C(1+|t|^2)^{\beta/2}$ for all $t\in\mathbb R$, $f$ is bounded and  $\int|\mathcal F f|$ is finite. If $F$ is not concave on $\mathbb{R}^+$, then there exists $C>0$ such that $T_n( F_{n}^{d}) \geq \sqrt nC$ with probability that tends to one as $n\to\infty$.
	\end{enumerate}
\end{theorem}

An interesting property of the test that can be seen in \eqref{convergence} is that the limiting distribution of our test statistic does not depend on the choice of $h_n$. Thus, the choice of this bandwidth can be made quite arbitrarily without altering the limit behavior of the test statistic under the null hypothesis.

\section{Calibration of the test}\label{sec: calib}

\subsection{Bootstrap calibration}\label{sec: bootstrap}

We will calibrate the test using bootstrap methods. For this task, we need to draw a bootstrap version of the observations in such a way  that the test statistic based on those observations mimics the behavior that the test statistic based on the original observations would have had if the null hypothesis were true. Hence, the bootstrap observations should be distributed as the sum between two independent variables, one of them having a concave distribution function that is close to that of $X$ under the null hypothesis, and the other one having the same distribution as $\varepsilon$. An initial suggestion of concave distribution function that could be used to that end, is the concave estimator $\mathcal MF_n^d$. 
However, this estimator is not necessarily a genuine distribution function on $\mathbb R^+$. We say that a function $H$ is a genuine distribution function on $\mathbb R^+$ if it is  a right-continuous non-decreasing function on $\mathbb R^+$ such that  $H(t)=0$ for any $t\le0$ and $\lim_{t\to\infty}H(t)=1$. To ensure that the  latter condition is fulfiled, we consider the estimator $\mathcal M\hat F_n^d$. Lemma \ref{lem: genuineDF} below establishes that this defines a genuine distribution function on $\mathbb R^+$. 
\begin{lemma}\label{lem: genuineDF}
	Let the assumptions of Theorem \ref{theo: ST} be satisfied. 
	Then,  
	$\mathcal M\hat F_n^d$ is a genuine distribution function
	on $\mathbb R^+$  with probability that tends to one. 
\end{lemma}
We turn to the construction of the bootstrap observations. Fix $m\ll n$ that tends to infinity as $n\to\infty$, and conditionally on the original sample $Y_1,\dots,Y_n$, consider a sample of $m$ iid random variables $X_1^*,\dots,X_m^*$ with common distribution function equal to $\mathcal M\hat F_n^d$, consider a sample $\varepsilon_1^*,\dots,\varepsilon_m^*$ from the same distribution as $\varepsilon$ and define $Y_i^*=X_i^*+\varepsilon_i^*$ for all $i$. Here, the sample $(\varepsilon_i^*)$ is independent of all other variables $(Y_1,\dots,Y_n,X_1^*,\dots,X_m^*)$. The reason why we choose $m\ll n$ is explained in Section \ref{sec: proofbootstrap}.

Now, compute the deconvolution estimator $F_m^{d*}$ based on observations $Y_1^*,\dots,Y_m^*$ and $T_m(F_m^{d*})$, the bootstrap version of the test statistic where  the functional $T_m$ maps a distribution function into the supremum distance between that distribution function and its least concave majorant on $\mathbb R^+$, times $\sqrt m$. Let $c^*_{\gamma,m}$ be the quantile of order $1-\gamma$ of $T_m(F_m^{d*})$ for some fixed prescribed level $\gamma\in(0,1)$. Note that $c^*_{\gamma,m}$ can be computed via Monte Carlo simulations. The critical region of the test is then
\begin{equation}\label{eq: critical}
	\left\{T_n(F_n^{d})>c_{\gamma,m}^*\right\}.
\end{equation}

To study the properties of the test, we need to explore the limit behavior of the bootstrap version of the test statistic, both under the null and the alternative hypotheses. This is the aim of the following theorem, which is a bootstrap version of Theorem \ref{theo: limit}. A key ingredient in the proof of Theorem \ref{theo: limit} is Theorem \ref{theo: ST}. Likewise, the proof of Theorem \ref{theo: bootstrap} relies on a bootstrap analogue of Theorem \ref{theo: ST}, namely Theorem \ref{theo: donsker}, which is stated in Section \ref{sec: generalization} and proved in Section \ref{sec: proofdonsker}. To prove that Assumptions {\bf B} of that theorem holds in our setting, we rely on results from the broader deconvolution literature, that requires stronger assumptions than for Theorem \ref{theo: limit}. Specifically, we assume the following conditions {\bf D'} and {\bf E'} instead of  {\bf D} and {\bf E}. \\

\noindent
{\bf Assumptions D'}
{\it 
	\begin{enumerate}
		\item  $f$  is  bounded, continuous and satisfies $\int
		\left\vert x\right\vert ^{8}f\left( x\right) dx<\infty$.
		
		\item $f\in H^{\alpha }\left( \mathbb{R}\right) $ for some $\alpha >1/2$ and \( \left| \left(\mathcal{F} f\right)^{(j)}(t) \right| \leq C(1 + |t|)^{ - j} \) for all \( 0 \leq j \leq 4,\ t\in\mathbb R\).
	\end{enumerate}
}

\noindent
{\bf Assumptions E'}
{\it 
	\begin{enumerate}
		\item$f_\varepsilon$  is  bounded and satisfies $\int
		\left\vert x\right\vert ^8f_\varepsilon\left( x\right) dx<\infty $.
		
		\item  $\varphi_{\varepsilon }$ has no vanishing
		points on the real line, and, for some \( \beta > 0 \) and $C>0$, the condition \eqref{eq: hypphiepsion} holds for all $t \in \mathbb{R}.$
	\end{enumerate}
}

We consider below a bootstrap version of the function $\tilde F_{n}^{d}$, see  \eqref{tildeFn}, that is involved in Theorem \ref{theo: ST}. Namely, we define
\begin{equation*}
	\widetilde F_{m}^{d*}\left( x\right) :=\int_{-\infty}^{x}\mathcal{F}^{-1}\left[ 
	\mathcal{F}K_{h_m}\frac{\varphi _{m}^*}{\mathcal{F}f_{\varepsilon }}\right] \left(
	t\right) dt
\end{equation*}
where  $\varphi_m^*(u)=m^{-1}\sum_{j=1}^me^{iuY_j^*}.$ Moreover, we
make  similar assumptions as in Theorem \ref{theo: ST} with $n$ replaced by $m$. However, 
in order to use results from \cite{hall2008estimation} on convergence of estimators of moments in the deconvolution model (see Section \ref{sec: momentdeconv}), we restrict attention to the specific kernel considered in that paper. Hence,
we now use a kernel \( K \) such that its Fourier transform is defined by
\begin{equation} \label{eq:kernel}
	\mathcal{F}K(t) = (1 - t^r)^s \cdot \mathds{1}_{\{|t| \leq 1\}},
\end{equation}
where \( r \geq 2 \) is an even integer and \( s \geq 1 \) is an integer. {In fact, we will assume that $r>4$ in order to ensure that  $\hat f_n$ has a bounded absolute moment of order $2+\delta$, see Section \ref{sec: fnhatmoments}. Moreover, in order to use Corollary \ref{cor}, we assume that  $h_n^{2+2\beta }n\rightarrow \infty $ as $n\to\infty$. It follows from Lemma \ref{lem: DE} below that the assumptions of Theorem \ref{theo: ST} are satisfied under the assumptions of the following bootstrap version of it.

	\begin{theorem}\label{theo: bootstrap}
		Suppose   that $m$ tends to infinity as $n\to\infty$ with $m\ll n$. 
		Assume conditions \textbf{D'} and \textbf{E'}  with $\beta<1/2$, $\alpha+3\gamma >2\beta +1$ for some $\gamma \in(\beta, 1/2)$. Assume that \( K \) is given by \eqref{eq:kernel} with $r>4$ and $rs>\lfloor\alpha+\gamma\rfloor$. Furthermore, let $h_{m}^{2\alpha+2\gamma }m\rightarrow 0$,  $h_{m}^{\rho }m\rightarrow \infty $ for some $\rho
		>4\beta -4\gamma +2$, $h_n^{2\alpha+2\gamma}n\to0$ and $h_n^{2+2\beta }n\rightarrow \infty $ as $n\to\infty$.
		Then,  we have:
		\begin{enumerate}
			\item
			Conditionally on $(Y_n)_{n\in\mathbb N}$ one has
			\begin{equation}\label{eq: cvboot}
				\sqrt m(\widetilde F_m^{d*}-\mathcal M\hat F_n^d) \rightsquigarrow \mathbb{G}_{F_0}\ \ in\ \ell_\infty(\mathbb{R}),
			\end{equation}%
			in probability as $n\to\infty$,
			where $F_0=\mathcal M F$ and $\mathbb{G}_{F_0}$ is as in Theorem \ref{theo: ST} with $F_0$ instead of $F$.
			\item If $F$ is concave  on $\mathbb R^+$ then conditionally on $(Y_n)_{n\in\mathbb N}$  one has
			\begin{equation*}
				T_m(F_m^{d*}) \stackrel{L}{\longrightarrow}\|\mathcal M_{F}'\mathbb{G}_{F}-\mathbb{G}_{F}\|_\infty
			\end{equation*}
			in probability 
			as $n\to\infty$.
			\item Whether or not $F$ is concave  on $\mathbb R^+$, the bootstrap test statistic $T_m(F_m^{d*})$ is bounded in probability conditionally on $(Y_n)_{n\in\mathbb N}$.
		\end{enumerate}
	\end{theorem}

	To understand the meaning of the conditional weak convergence in probability in \eqref{eq: cvboot}, one can consider the Prohorov distance, that metrizes weak convergence of probability measures on separable spaces, see e.g. \cite[page 72]{billingsley2013convergence}. The Prohorov distance $\pi(\mu,\nu)$ between two probability measures $\mu$ and $\nu$ on $\ell_\infty(\mathbb{R}^+)$ (the separable space that we consider in Theorem \ref{theo: bootstrap} above),  is defined as the infimum of those $\epsilon>0$ for which the two inequalities
	\begin{equation*}
		\mu(A)\leq\nu(A^\epsilon)+\epsilon,\ \nu(A)\leq\mu(A^\epsilon)+\epsilon
	\end{equation*}
	hold for all Borel sets $A$, where
	\begin{equation*}
		A^\epsilon=\left\{f\in\ell_\infty(\mathbb{R}^+),\ \inf_{g\in A}\|f-g\|_\infty<\epsilon\right\}.
	\end{equation*}
	This distance metrizes weak convergence in the sense that a sequence of probability measures $(\mu_n)_{n\in\mathbb{N}}$ weakly converges to a probability measure $\mu$ on $\ell_\infty(\mathbb{R}^+)$ if and only if $\pi(\mu_n,\mu)\to 0$ as $n\to\infty$. If we denote by $\mu_n^*$ the conditionnal distribution of 
	$\sqrt m(F_m^{d*}-\hat{F}_n^d)
	$ 
	and by $\mu$ the distribution of $\mathbb{G_F}$, then the meaning of \eqref{eq: cvboot} is that \begin{equation}\label{eq: prohorov}
		\pi(\mu_n^*,\mu)  \stackrel{\mathbb{P}}{\longrightarrow}0,\ as\ n\to\infty.
	\end{equation}
	
	We now consider consistency of the test under both the null and the alternative hypotheses. Hence, we assume below that the assumptions of Theorem \ref{theo: bootstrap} are satisfied. Thanks to Lemma \ref{lem: DE} below, we know that, in that case, the assumptions of Theorems \ref{theo: ST} and \ref{theo: limit} are satisfied as well.
	
	\subsection{Asymptotic properties of the test}
	We assume without loss of generality that $\gamma<1/2$ since levels of tests are always choosen with this constraint (recall that a typical level is $\gamma=0.05$). The following theorem, which is proved in Section \ref{sec: prooflevel}, shows that the test has the prescribed asymptotic level $\gamma$ at each point $F$ in the null hypothesis that has at least an affine part. The proof of the theorem uses the fact that, with probability one, the process  $\mathbb G_F$ cannot be affine on a given non-degenerate interval on which $F$ is affine. We show in Section \ref{sec: prooflevel} that this property holds in the case where the distribution of $\epsilon$ is supported on the whole real line, hence we make this assumption in the theorem. We conjecture that this holds under more general conditions, but due to the very complex form of the covariance function of $\mathbb G_F$, proving it is out of the scope of the paper.  We can notice that in the case of zero noise, the limit process is $B\circ F$ where $B$ is standard Brownian bridge, and it is well known that with probability one, $B$ (hence also $B\circ F$) cannot be affine on a given interval on which $F$ is affine. Hence, the conjecture generalizes to the case of measurement error a property that is well know under no error.
	\begin{theorem}\label{lem: level}
		Assume $\gamma<1/2$, the distribution of $\epsilon$ is supported on $\mathbb R$, and the assumptions of Theorem \ref{theo: bootstrap} are satisfied. If $F$ is concave  on $\mathbb R^+$ with at least an affine part, then 
		\begin{equation}\label{eq: level}
			\lim_{n\to\infty}\mathbb P\left(T_n(F_n^{d})>c_{\gamma,m}^*\right)=\gamma.
		\end{equation}
	\end{theorem}

	On the other hand, if $F$ is strictly concave and nowhere affine, the limit of the test statistic is zero. In particular, the limit distribution is not continuous and we doubt that \eqref{eq: level} holds. In fact, our conjecture 
	is that the limit is equal to zero.  
	The conjecture is based on the following heuristics. Since the original test statistic $T_n(F_n^d)$ converges to zero in probability, it is expected that there exists a sequence of positive numbers that tend to infinity and such that $v_nT_n(F_n^d)$ is stochastically bounded and bounded away from zero with hight probability. In such a case, one can expect that the bootstrap test statistic has a similar behaviour whence $v_mT_m(F_n^{d*})$ is bounded from above and away from zero with hight probability. Then, its $(1-\gamma)$-quantile is of order of magnitude $v_m^{-1}$. Since $m\ll n$ one has $v_m^{-1}\gg v_n^{-1}$ whence the probability that $T_n(F_n^d)$, which is of order $v_n^{-1}$,  exceeds the quantile of order $v_m^{-1}$ tends to zero 	as $n$ tends to infinity.
	If these intuitions are correct, we may expect that the probability of rejection tends to 0, whence the test has asymptotic prescribed level $\gamma$ at any point $F$ in the null hypothesis, and is conservative in case of strict concavity. This is  confirmed by our simulation analysis in Section \ref{sec: simul}. 
	
	
	}

	Assume now that  $F$ is not concave. For arbitrary $C>0$ one has
	\begin{eqnarray*}
		\mathbb P\left(T_n(F_n^{d})>c_{\gamma,m}^*\right)&\geq&
		\mathbb P\left(c_{\gamma,m}^*\leq\sqrt n C\right) - \mathbb P\left(T_n(F_n^{d})\leq \sqrt n C\right). 
	\end{eqnarray*}
	Since \eqref{eq: lapfini} below holds, we can apply Theorem \ref{theo: limit}, so one can choose $C>0$ such that the second probability on the right hand side tends to zero as $n\to\infty$. On the other hand, it follows from Theorem \ref{theo: bootstrap} that $T_m(F_m^{d*})$ is bounded in probability conditionally, which implies that $c_{\gamma,m}^*$ is bounded in probability, conditionally and also unconditionally. Hence, the first term on the right hand side  converges to one as $n\to\infty$,
	\begin{eqnarray*}
		\lim_{n\to\infty}\mathbb P\left(T_n(F_n^{d})>c_{\gamma,m}^*\right)=1.
	\end{eqnarray*}
	In other words, the test defined by the critical region in \eqref{eq: critical} is consistent in the sense that it has prescribed asymptotic level, and its asymptotic power is equal to one at any point of the alternative.

	\subsection{About least favorable hypothesis}\label{sec: lfh}
	A common way to calibrate tests without bootstrap is to exhibit a least favorable null distribution $F_0$, i.e., a function satisfying the null hypothesis such that, under the null, $T_n(F_n^d)$ is asymptotically stochastically smaller than $T_{n0}$ where $T_{n0}$ denotes the statistic obtained when the true underlying function equals $F_0$. Since $F_0$ is known, one can easily simulate $T_{n0}$ and compute its $1-\gamma$ quantile using Monte Carlo. Then, the test that rejects the null hypothesis if $T_n(F_n^d)$ exceeds this quantile has (asymptotic) level $\gamma$, but it has smaller power than that  of the test calibrated using bootstrap. Anyway, this kind of calibration is appealing since it is easier to implement than the bootstrap, so  it would be interesting to see whether or not it could be  implemented  in our case. 

	
	When testing concavity of a distribution function using the $\ell_\infty$-distance between the empirical distribution function and its least concave majorant as a test statistic, based on an i.i.d. sample from that distribution, the uniform  is known to be least favorable, see  \cite{kulikov2004testing}. Hence, it is natural to wonder if uniform distributions are least favorable in our context. Unfortunately, it seems that it is not the case. 
	Indeed, uniform distributions just differ in terms of scale parameters. In the classic framework, where there is no measurement error, scaling does not make any difference, so one can conveniently choose the uniform on $[0,1]$. However, in our setting this is not possible. If we rescale $X$, this corresponds to rescaling $\epsilon$, that is, changing the noise-to-signal ratio, which has a crucial impact on the estimation accuracy, and consequently, on the power of the test. Since the noise-to-signal ratio is generally unknown, it is not possible to know which uniform distribution is the least favourable in a given case. Moreover, even assuming that the variance of $X$ is known, we are not able to prove that affine functions are least favorable because the way the distribution of our estimator depends on our function $F$ is very intricate. Indeed, the limit covariance $\Sigma_{s,t}$ depends on $f_Y=f_\varepsilon\star f$ which cannot be expressed in terms of $f_\varepsilon\star f_0$, with $f_0$ indicating the density of the uniform distribution on the support of $f$. 
	

	A simple way to circumvent the absence of a least favorable case,
	at the cost of conservativeness, is to reject the null hypothesis
	whenever $T_n(F_n^d)$ exceeds a deterministic threshold $v_n$,
	where $v_n \to \infty$ and $v_n = o(\sqrt{n})$ as $n \to \infty$.
	For instance, one may take $v_n = \log n$. Under the null hypothesis, Theorem~\ref{theo: limit} establishes the limit distribution of the test statistic, which is almost surely finite, and consequently
	$T_n(F_n^d)=O_P(1)$.
	Then, $P(T_n(F_n^d) > v_n) \to 0$, hence, the resulting test has asymptotic level equal to zero and is conservative. Under any fixed alternative,
	the test statistic diverges to infinity at rate $\sqrt{n}$.
	Since $v_n = o(\sqrt{n})$, it follows that the test is consistent and has asymptotic power equal to one at every point in the alternative.

	\section{Simulations}\label{sec: simul} 
	
	In this section,  we implement the estimator, evaluate its performance, and assess the test by numerically calculating its power for different target distributions. To construct $f^d_{n}$, we employ a kernel with the following Fourier transform 
	\[
	\mathcal{F}K(t) = (1 - t^{6})\,\mathds{1}_{[-1,1]}(t),
	\]  
	which, according to Lemma \ref{lem: DE}, satisfies Assumption \textbf{K} with $L = 5$.  
	The bandwidth $h_n$ is chosen via the plug-in method of \cite{delaigle2004practical}, which estimates an asymptotic approximation of the mean integrated square error to obtain a data-driven optimal bandwidth.
	The remainder of this section is organized as follows. In Subsection \ref{sim:mse}, we evaluate the performance of the proposed estimators from Section \ref{sec: concaveestim}. In Subsection \ref{sim:test} we compute and examine the rejection rates of our test from Section \ref{sec: calib} under different scenarios.
	All computations are implemented in \texttt{R}; the code is available at
	{https://github.com/MohammedEssalih/Concave-CDF-under-measurement-error}.
	


	\subsection{Performance of $\mathcal{M}\hat{F}_n^d$ in $\mathcal{H}_0$}\label{sim:mse}

	This section evaluates the performance of the constrained estimator $\mathcal{M}\hat{F}_n^d$ in terms of its mean squared error (MSE). The objective is to determine whether the constrained estimator $\mathcal{M}\hat{F}_n^d$ provides improvements over its unconstrained counterpart $\hat{F}_n^d$, when $F$ is concave. To this end, we analyze the ratio of the MSE of $\mathcal{M}\hat{F}_n^d$ and $\hat{F}_n^d$ computed at different quantiles, based on $500$ simulation runs. Here, $\mathcal{M}\hat{F}_n^d$ outperforms $\hat{F}_n^d$ when the MSE ratio is smaller than one.
	
	Denote with $W(a,b)$ a Weibull distribution with shape and scale parameters $a>0$ and $b>0$, respectively. The distribution function of the $W(a,b)$ is concave when $a\leq1$ and not concave for $a>1$ (for $a>1$ the function is convex to the left of an inflection point and concave to the right). We also consider the Beta distribution with shape parameters $a,b>0$, denoted with $B(a,b).$ Setting $b=1$ for simplicity, the distribution function of $B(a,1)$ is concave for $a\leq 1$, convex for $a\geq 1$, and the case $a=1$ corresponds to the uniform (the distribution function is the identity on $[0,1]$). 
	
	
	To ensure a comprehensive evaluation, the estimation is performed over quantile levels within the interval $(0, 1)$. The bounds $0$ and $1$ are omitted because the MSE of both estimators are equal to zero at those points (by definition of $\hat F_n^d$ and Lemma \ref{lem: genuineDF}), whence the ratio is undefined there.
	
	For simplicity, we restrict our attention to the case where the measurement error follows a Laplace distribution, with mean $0$ and standard deviation $\sigma_{\varepsilon}$. The standard deviation is chosen to achieve a specific noise-to-signal ratio (NSR), following \cite[p.~483]{dattner2013estimation}, defined as $$\text{NSR} = \frac{\sigma_{\varepsilon}}{\sigma_X},$$
	which quantifies the degree of error contamination. Three levels of contamination, corresponding to NSR values $0.1$, $0.2$, and $0.5$, are examined across sample sizes $n \in \{20, 50, 100\}$.
	
\begin{figure}[ht]
		\centering
		\begin{subfigure}[b]{\textwidth}
			\centering
			\includegraphics[width=4.5 cm,height=4 cm]{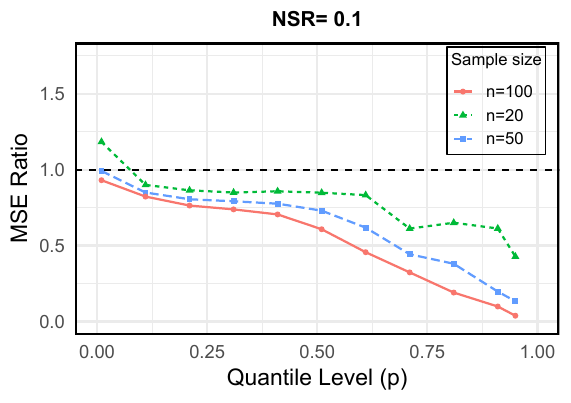}
			\includegraphics[width=4.5 cm,height=4 cm]{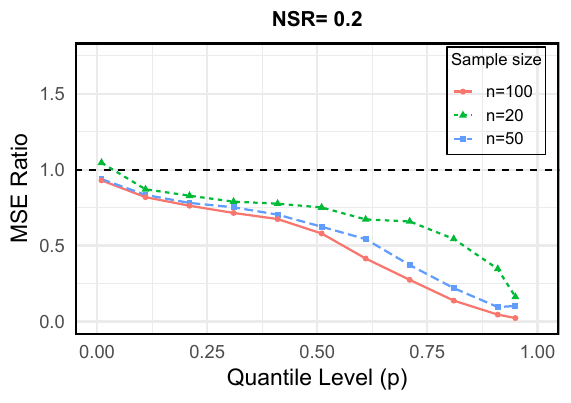}
			\includegraphics[width=4.5 cm,height=4 cm]{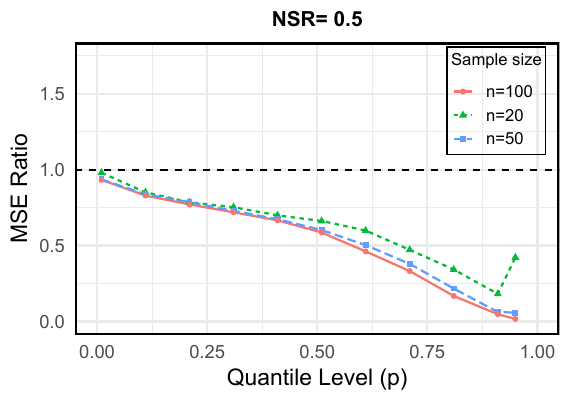}
			\caption{ $W(0.75,1)$.}
			\label{MSE_Weibull}
		\end{subfigure}
		
		\vspace{0.5em}
		
		\begin{subfigure}[b]{\textwidth}
			\centering
			\includegraphics[width=4.5 cm,height=4 cm]{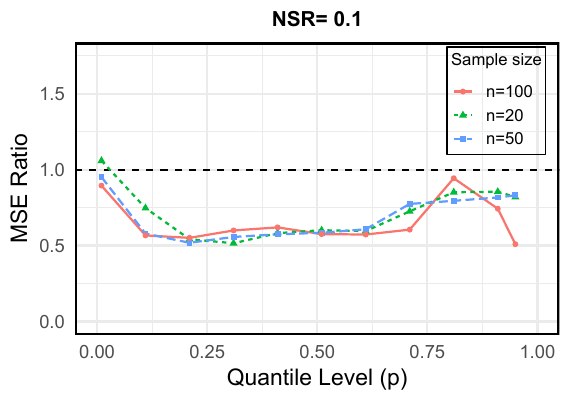}
			\includegraphics[width=4.5 cm,height=4 cm]{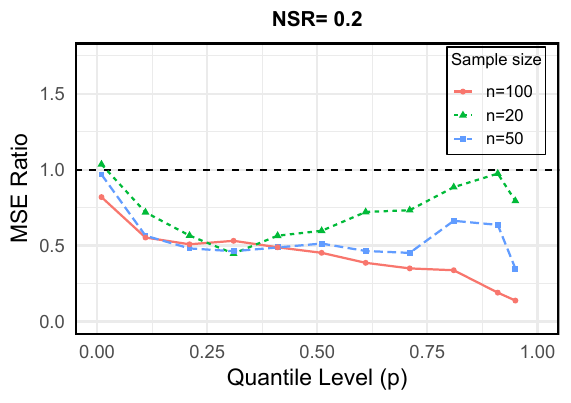}
			\includegraphics[width=4.5 cm,height=4 cm]{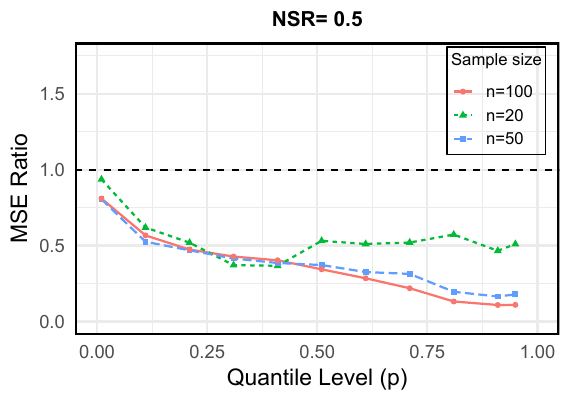}
			\caption{$B(0.75, 1)$}
			\label{subfig:MSE_beta}
		\end{subfigure}
		
		\caption{Simulation results. Mean squared error (MSE) for the Beta $B(0.75,1)$ and Weibull $W(0.75,1)$ models, across NSR levels of $10\%$, $20\%$, and $50\%$ (from left to right) and sample sizes $n = 20, 50,$ and $100$.}
		\label{MSE}
	\end{figure}

	The obtained results for the Weibull distribution $W(0.75, 1)$ and the Beta distribution $B(0.75, 1)$ are summarized in Figure~\ref{MSE}.

	The results presented in Figure~\ref{MSE} show that the constrained estimator $\mathcal{M}\hat{F}_n^d$ can substantially outperform its unconstrained counterpart, especially in the right tail. 

	\subsection{Rejection rates for $T(F_n^d)$}\label{sim:test}
	
	In this section, we proceed as follows. First, we describe the error distribution used in the test. Next, we implement the test in the case where the distribution function is strongly concave and nowhere affine, using the Beta and Weibull distributions. We also consider the weakly concave case, where the distribution function contains affine segments. Finally, we present an example illustrating how the presence of error can alter the shape of the distribution function, and we apply the test to examine whether it can detect this misleading effect.
	
	\subsubsection{The distribution of the error used for the test} 
	The assumptions of Theorem \ref{theo: ST} require a smoothness condition for the error distribution. Such a condition must hold for some $\beta < 1/2$, excluding, for instance, the Laplace distribution, which corresponds to $\beta = 2$.  
	
	Differently, the symmetric-gamma distribution $\text{SG}(\beta, \theta)$, with shape parameter~$\beta>~0$, scale parameter~$\theta>0$, and characteristic function  $ (1 + \theta t^2)^{-\beta}$, $t\in \mathbb{R}$, satisfies the assumptions for $\beta \in (0, 1/4)$, see Example~3 of \cite{sohl2012uniform}. Thereore, we assume that the error distributiuon is a mixture between Laplace and symmetric-gamma, with density
	\[
	f_{\varepsilon}(x) \;=\; p \, f_{\text{SG}}\!\left(\beta,\,\theta\right) 
	\;+\; (1-p) \, f_{\text{Lap}}\!\left(0,\,\sigma_{\text{Lap}}\right),\quad p\in(0,1),
	\]
	where $f_{\text{SG}}$ denotes the density of a symmetric-gamma distribution, and $f_{\text{Lap}}$ denotes the density of a Laplace distribution with mean~$0$ and scale parameter~$\sigma_{\text{Lap}}>0$. It is easy to see that this mixture model satisfies the assumptions because of the symmetric-gamma component. 
	Moreover, the Laplace component allows altering the shape of the convolution distribution function in a more effective way compared to the symmetric-gamma alone, providing more interesting cases of study. Indeed, the distribution function of the symmetric-gamma distribution is strictly concave in $(0,\infty)$, when convolved with some $F$, the resulting distribution is typically still concave regardless of the shape of $F$.
	
	
	On the one hand, the choice of the symmetric-gamma parameter is not particularly critical: if $p$ is chosen very small, the overall shape of the error distribution is mainly determined by the second component of the mixture. On the other hand, the smoothness of the resulting mixture remains comparable to that of the symmetric-gamma distribution whenever the second component in the mixture is smoother than the symmetric-gamma. 

	To compute the NSR, the variance of~$\varepsilon$ is $	\operatorname{Var}(\varepsilon)
	= p \cdot \beta(1+\beta)\theta^{2} + (1-p) \cdot 2\sigma_{\text{Lap}}^{2}$. Setting \(\theta = 0.25\), \(\beta = 0.24\), and \(p = 0.01\), $\sigma_\varepsilon$ just depends on $\sigma_{\text{Lap}}$. To control the NSR, we choose \(\sigma_{\text{Lap}}\) so that	$\sigma_\varepsilon \;=\; \mathrm{NSR} \times \sigma_X,$
	where \(\sigma_X \) is the standard deviation of $X$, and \(\mathrm{NSR} \in \{0.1,\,0.2,\,0.5\}\).

	\subsubsection{The power of the test}
	
\begin{figure}[ht]
		\centering
		\begin{subfigure}[b]{\textwidth}
			\centering
			\includegraphics[width=0.32\textwidth]{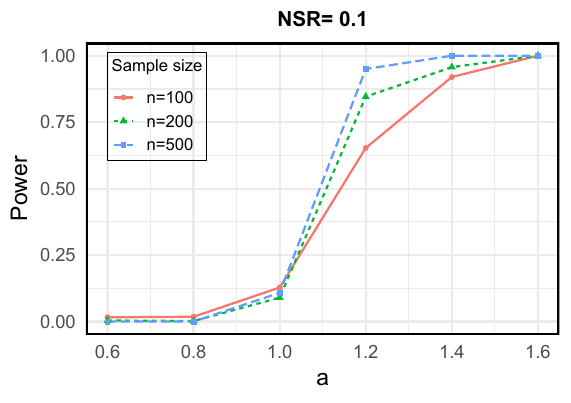}
			\includegraphics[width=0.32\textwidth]{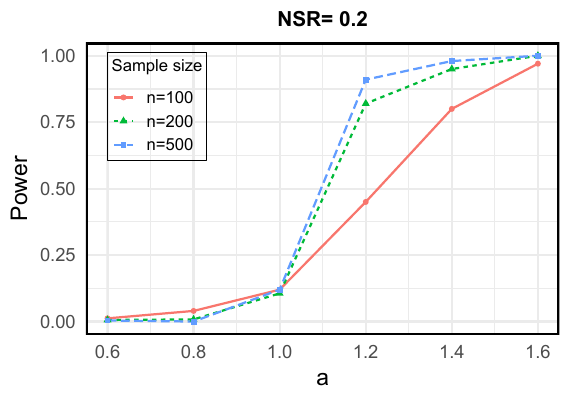}
			\includegraphics[width=0.32\textwidth]{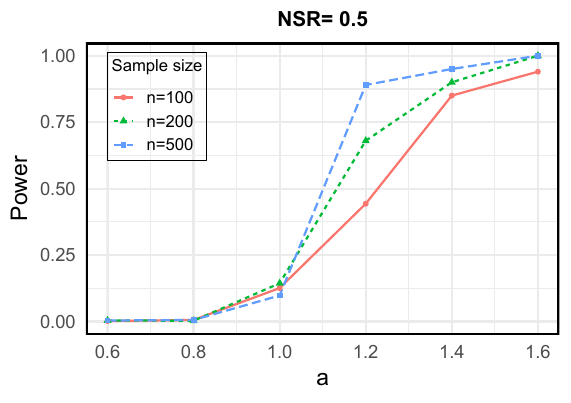}
			\caption{$\mathrm{Beta}(a, 1)$ with $a = 0.6, 0.8, \ldots, 1.6$.}
			\label{subfig:beta}
			\vspace{0.5em}
			\begin{subfigure}[b]{\textwidth}
				\centering
				\includegraphics[width=0.32\textwidth]{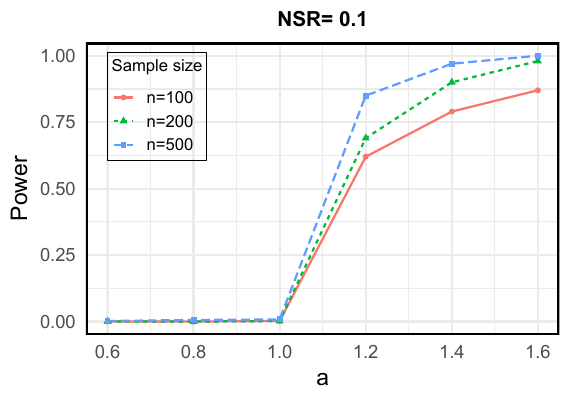}
				\includegraphics[width=0.32\textwidth]{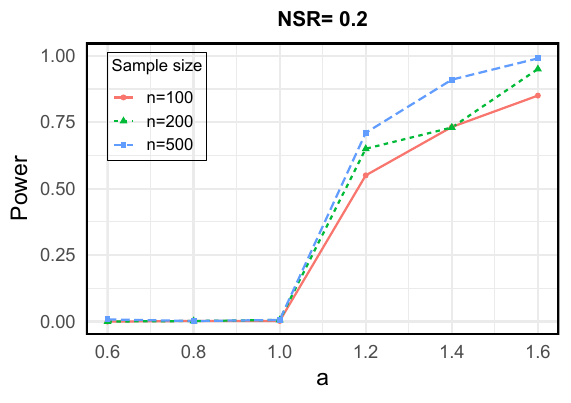}
				\includegraphics[width=0.32\textwidth]{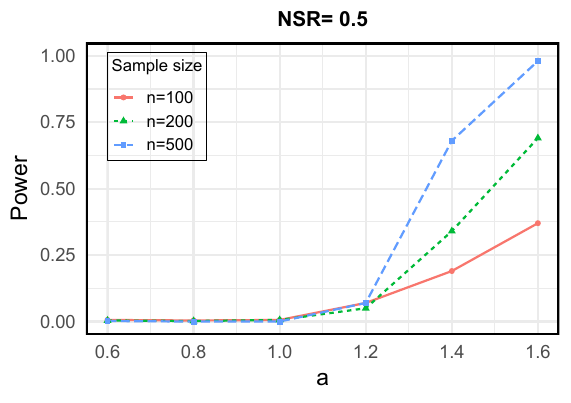}
				\caption{Weibull $W(a, 1)$ with $a = 0.6, 0.8, \ldots, 1.6$.}
				\label{subfig:weibull}
			\end{subfigure}

		\end{subfigure}
		
		\caption{Rejection rates for $\mathrm{Beta}(a, 1)$ and $W(a, 1)$ distributions for different values of $a$, NSR, and $n$.}
		\label{fig:power}
	\end{figure}

	The test performance is evaluated using $M=500$ Monte Carlo replications, $m=\lfloor n^{0.9}\rfloor$, plug-in-selected bandwidth $h_m$, $B=300$ bootstrap samples, and significance level $\gamma=0.1$. First, we consider \( X \sim B(a, 1) \) with the shape parameter \( a \) ranging from 0.6 to 1.6. The results are plotted and summarized in Figure~\ref{subfig:beta}, with the x-axis representing $a$. Recall that the distribution function of the \( B(a, b) \) is concave for \( a \leq 1 \) and non-concave for \( a > 1 \), with \( B(1,1) \) corresponding to the uniform. 
	
	Figure~\ref{subfig:beta} indicates that, under $\mathcal{H}_0$ (i.e., $a \leq 1$), the rejection rates of $T(F_n^d)$ remain below the nominal level $\gamma = 0.1$. In particular, they tend to 0 in the strictly concave case ($a<1$) and to $\gamma=0.1$ in the affine case ($a=1$), which supports our heuristics in the paragraph after \eqref{eq: level}. Moreover, the empirical power tends to $1$ under \( \mathcal{H}_1 \) as the sample size \( n \) increases. As expected, for larger (smaller) values of the NSR, the performance of the test seems to deteriorate (improve). This is quite logical, as it becomes harder to detect which hypothesis is true when the measurement error gets larger.
	We also examine the Weibull distribution \( W(a, b) \), which belongs to \( \mathcal{H}_0 \) for \( a \leq 1 \) and to \( \mathcal{H}_1 \) otherwise. 
	Figure~\ref{subfig:weibull} shows a similar pattern, with the distinction that this cumulative distribution has no affine segments. Here, the power of the test approaches $0$ for $a \leq 1$ and $1$ for $a > 1$.

	We now consider a second scenario within $\mathcal{H}_0$, namely, a distribution function which is affine on a non-empty interval and strictly concave elsewhere. To illustrate this, we consider the following mixture $
	X \sim 1/2 \cdot \mathrm{U}[0,1] \;+\; 1/2 \cdot \bigl(1 + \mathrm{Exp}(1)\bigr),$
	where \(\mathrm{U}[0,1]\) denotes the uniform distribution on $[0,1]$ \(\mathrm{U}[0,1]\) and $\mathrm{Exp}(1)$ is the standard exponential 
	distribution with scale parameter equal to 1. This mixture distribution is affine on \([0,1]\) and strictly concave on \((1, \infty)\). The simulation results, reported in Table~\ref{Weak}, show that the rejection rates converge to the significance level $\gamma$ as $n$ increases, which also supports our heuristics in the paragraph after \eqref{eq: level}.

	\begin{table}[h!]
		\centering
		\begin{tabular}{|c|c|c|c|}
			\hline
			$n$ & NSR = 0.1 & NSR = 0.2 & NSR = 0.5 \\
			\hline
			100 & 0.115 & 0.085 & 0.060 \\ \hline
			200 & 0.103 & 0.095 & 0.080 \\ \hline
			500 & 0.1005 & 0.0995 & 0.098 \\
			\hline
		\end{tabular}
		\caption{Rejection rates for the mixture $X \sim 0.5 \cdot \mathrm{U}[0,1] \;+\; 0.5 \cdot (1 + \mathrm{Exp}(1))$ distribution for different sample sizes $n$ and NSR values.}
		\label{Weak}
	\end{table}

	\subsection{The performance of the test when the error alters the shape}
	In some situations, the presence of measurement error affects the shape of the convolution substantially. For example, when $X$ follows the mixture distribution $X \sim 0.2\,W(3,1) + 0.8\,B(0.5,0.75),$ then $F$ satisfies $\mathcal H_1$. However, after convolution with measurement error, the distribution function of $X+\varepsilon$ may fall into $\mathcal H_0$, especially when the NSR is large. In such cases, it is essential to apply the proposed test in order to discover the shape of the true underlying distribution $F$.
	
	The performance of the test is reported in Table~\ref{T3} for different NSR levels and sample sizes~$n$. The results show that the NRS has a big impact on the performance. However, as $n$ increases, the test is able to detect the non-concavity of $F$, especially for moderate NRS values. In situations like this, our test is especially useful, as it helps discovering the true shape of a distribution function, removing the impact of measurement error.

	\begin{table}[h!]
		\centering
		\begin{tabular}{|c|c|c|c|c|c|c|}
			\hline
			\textbf{$NSR$} & \textbf{$n=20$} & \textbf{$n=50$} & \textbf{$n=100$} & \textbf{$n=200$} & \textbf{$n=500$} & \textbf{$n=1000$} \\ \hline
			$0.1$ & 0.012 & 0.063 & 0.109 & 0.211 & 0.594 & 0.781 \\ \hline
			$0.2$ & 0.018 & 0.051 & 0.082 & 0.147 & 0.358 & 0.601 \\ \hline
			$0.5$ & 0.026 & 0.039 & 0.045 & 0.072 & 0.101 & 0.176 \\ \hline
		\end{tabular}
		
		\caption{Rejection rates for the mixture $X \sim 0.2\,W(3,1) + 0.8\,B(0.5,0.75)$ distribution for different sample sizes $n$ and $NSR$ values.}
		\label{T3}
	\end{table}

	\appendix
	
	\section{Generalization of the Donsker theorem in  \cite{sohl2012uniform} to a bootstrap setting}
	\label{sec: generalization}
	Recall that in the setting of  \cite{sohl2012uniform}, one observes independent and identically distributed random variables  $Y_j=X_j+\varepsilon_j$ for $j=1,\dots,n$ where $X_j$ and $\varepsilon_j$ are independent of each other, the distribution of $\varepsilon_j$ is supposed to be known and the aim is statistical inference on the distribution of $X_j$. In this deconvolution model, the random variables $X_1,\dots,X_n$ are not observed, but only $Y_1,\dots,Y_n$, which contains an additive error. In this model, kernel estimators $\widehat{ \vartheta_t}$ are considered for linear functionals
	\begin{eqnarray*}
		t\mapsto \mathcal \vartheta_t:= \int\zeta(x-t)f_X(x)dx,
	\end{eqnarray*}
	with $f_X$ being the probability density function of the $X_i$s, 
	where the special case $\zeta:=\mathds{1}_{(-\infty,0]}$ leads to the estimation of the distribution function of the $X_i$s. Theorem 1 in  \cite{sohl2012uniform} states that $\sqrt n(\widehat{ \vartheta_t}-\vartheta_t)_{t\in\mathbb{R}}$ converges in law to a centered Gaussian process $\mathbb{G}$ in $\ell_\infty(\mathbb{R})$, whose covariance depends on the functional $\zeta$ and on the distribution of the errors $\varepsilon_j$. Here,
	\begin{equation*}
		\widehat{ \vartheta_t}=\int\zeta(x-t)\mathcal F^{-1}\left[\mathcal FK_{h_n}\frac{\varphi_n}{\varphi_\epsilon}\right](x)dx,
	\end{equation*}
	where $K_{h_n}$ and $\varphi_n$ are defined as in \eqref{fn}, and $\varphi_\varepsilon$ denotes the Fourier transform of the density of $\varepsilon$.
	
	The aim of the current section is to establish a bootstrap version of that theorem. 
	Hence, we assume that we have at hand an estimator $\hat f_n$ for the density of the $X_j$s based on observations $Y_1,\dots,Y_n$. With $(m)$ a sequence of positive integers that may depend on $n$, conditionally on $Y_1,\dots,Y_n$, we build independent and identically distributed variables $X_1^*,\dots,X_m^*$ with density function $\hat f_n$. (Note that we assume here that $\hat f_n$ is a genuine density function, which means that it is non-negative and integrates to one). Independently of all those variables, we build independent and identically distributed variables $\varepsilon_1^*,\dots,\varepsilon_m^*$ with the same distribution as the original errors $\varepsilon_j$s. Then, conditionally on the original observations $Y_1,\dots,Y_n$, we obtain  independent and identically distributed bootstrap observations
	\begin{equation*}
		Y_j^*=X_j^*+\varepsilon_j^*,\quad j=1,\dots,m.
	\end{equation*}
	We then define the estimator in the bootstrap world as follows 
	\begin{equation*}
		\widehat{ \vartheta_t^*}=\int\zeta(x-t)\mathcal F^{-1}\left[\mathcal FK_{h_m}\frac{\varphi_m^*}{\varphi_\epsilon}\right](x)dx,
	\end{equation*}
	where $\varphi_m^*(u)=m^{-1}\sum_{j=1}^me^{iuY_j^*},$ $u\in\mathbb{R}$ and $\zeta:\mathbb R\to\mathbb R$ is a given function that belongs to $Z^{\gamma_s,\gamma_c}$ for some $\gamma_s,\gamma_c>0$, where (as in \cite{sohl2012uniform}) $Z^{\gamma_s,\gamma_c}$ is the set of functions $\zeta$ such that $\zeta=\zeta^c+\zeta^s$ for some $\zeta^s\in H^{\gamma_s}(\mathbb R)$ that is compactly supported as well as $\langle x\rangle^\tau(\zeta^c(x)-a(x))\in H^{\gamma_c}(\mathbb R)$ for some $\tau>0$ and some $a\in C^\infty(\mathbb R)$ such that $a'$ is compactly supported.
	Similar to Assumption {\bf D} in Section \ref{sec: assumptions} (which is Assumption 1 in \cite{sohl2012uniform}) we make the following assumptions for the bootstrap version. \\

	\noindent
	{\bf Assumption B}
	{\it 
		\begin{enumerate}
			\item
			The estimator $\hat f_n$ is a genuine density function such that
			\begin{equation}\label{eq: fnhatbounded}
				\sup_x\hat f_n(x)=O_P(1)
			\end{equation}
			and 
			\begin{equation}\label{eq: fnhatmoments}
				\int |x|^{2+\delta}\hat f_n(x)dx=O_P(1)
			\end{equation}
			for some arbitrarily small $\delta>0$. 
			\item 
			One has
			\begin{equation}\label{eq: assumption fii}
				\int\langle u\rangle^{2\alpha}|\mathcal F \hat f_n(u)|^2 du=O_p(1)
			\end{equation}
			for some $\alpha>0$, where $\langle u\rangle=(1+u^2)^{1/2}$.
			\item $\Vert \hat f_n-f\Vert_\infty$ converges to zero in probability as $n\to\infty$. 
		\end{enumerate}
	}

	Under these assumptions, we have the following Donsker theorem, which generalizes Theorem 1 in \cite{sohl2012uniform} to a general bootstrap setting.
	
	\begin{theorem}\label{theo: donsker}
		Suppose that $m$ tends to infinity as $n\to\infty$.
		Assume  conditions \textbf{E} and \textbf{B} as well as $\zeta\in Z^{\gamma_s,\gamma_c}$ with $\gamma_s>\beta$, $\gamma_c>(1/2\vee\alpha)+\gamma_s$ and $\alpha+3\gamma_s >2\beta +1$. Furthermore, let the kernel $K$ satisfy conditions \textbf{K} with $L = \left\lfloor \alpha + \gamma_s \right\rfloor$. Let $h_{m}^{2\alpha+2\gamma_s }m\rightarrow 0$ and if $\gamma_s\leq\beta+1/2$ let in addition $h_{m}^{\rho }m\rightarrow \infty $ for some $\rho
		>4\beta -4\gamma_s +2$. Then,  contionally on $(Y_n)_{n\in\mathbb N}$ one has
		\begin{equation*}
			\sqrt m(\widehat{\vartheta_t^*}-\vartheta_t^*)_{t\in\mathbb{R}} \rightsquigarrow\mathbb{G}_F\ \ in \ \ell_\infty(\mathbb {R})
		\end{equation*}
		in probability 
		as $n\to\infty$, where $\vartheta_t^*=\int\zeta(x-t)\hat f_n(x)dx$ and $\mathbb{G}_F$ is the same centered Gaussian Borel random element  in $\ell_\infty(\mathbb{R})$ as in Theorem \ref{theo: ST}.
	\end{theorem}
	
	\section{Proofs for Sections \ref{sec: context} and \ref{sec: concaveestim}}\label{sec: proofcvunif}
	
	\subsection{Proof of Lemma \ref{lem: welldef}}
	{
		We write $h$ insted of $h_n$.
		The proof relies on the results from \cite{Liu1989} concerning the strong consistency of  density estimators. Since $\mathcal F K$ is supported on $[-1,1]$,  the estimator in (\ref{fn}) aligns with the truncated density estimator in \cite{Liu1989}. Specifically, if the Fourier transform $\mathcal{F}K$ is compactly supported in $\left[ -1,1 \right]$, with $K_h(x) = h^{-1} K(x/h)$, then $\mathcal{F}K_h(t) = \mathcal{F}K(th)$ is also compactly supported, with support in $\left[ -1/h,1/h \right]$. Consequently, the estimator in (\ref{fn}) can be expressed with a finite integral instead of over $(-\infty, +\infty)$:
		\begin{equation*}
			f_{n}^{d}\left( x\right) =\frac{1}{2\pi }\int_{-M}^{M}\exp \left(- itx\right) 
			\mathcal{F}K\left( th\right) \frac{\varphi _{n}\left( t\right) }{\mathcal{F}
				f_{\varepsilon }\left( t\right) }dt\text{ where }M=1/h.
		\end{equation*}
		Arguing as in the proof of Theorem 3.1 in  \cite{Liu1989} with $M=1/h$ proves that 
		\begin{equation*}
			\mathbb E\int_{-M}^M|\mathcal F f_n^d(t)-\mathcal F f(t)|dt\leq \frac{M\sqrt 2}{\pi\sqrt n\min_{|t|\leq M}|\varphi_\varepsilon(t)|}+o(1),
		\end{equation*}
		and that 
		\begin{equation*}
			\mathbb E\|f_n^d-f\|_\infty\leq \frac{1}{2\pi}\mathbb E\int_{-M}^M|\mathcal F f_n^d(t)-\mathcal F f(t)|dt+ \int_{|t|>M} |\mathcal{F}f(t)| \, dt. 
		\end{equation*}
		Under our assumptions one has $1/|\varphi_\varepsilon(t)|\leq CM^{\beta}$ for some positive constant $C$ and all $t\in[-M,M]$ so we have
		\begin{equation*}
			\frac{M}{\sqrt n\min_{|t|\leq M}|\varphi_\varepsilon(t)|}\leq C\frac{M^{1+\beta}}{\sqrt n}
		\end{equation*}
		where the right-hand side converges to zero under the assumption that $h^{ 2+2\beta }n\rightarrow \infty $.
		Moreover,  $\int |\mathcal{F}f(t)| \, dt $ is finite by assumption, so $ \int_{|t|>M} |\mathcal{F}f(t)| \, dt $ converges to zero as $n\to\infty$ and we get
		\begin{equation*}
			\mathbb E\| f_n^d-f\|_\infty=o(1).
		\end{equation*}
		It then follows from Markov's inequality that
		\begin{equation*}
			\| f_n^d-f\|_\infty=o_p(1).
		\end{equation*}
		Moreover, for all $t\geq 0$,  one has 
		$$
		|F^d_{n}(t) -F(t) | \leq \int_{0}^{t}|f^d_n(u) -f(u) |du \leq t\|f^d_{n} -f \|_\infty
		,	$$
		where the supremum on the right-hand side tends to zero.
		Then  \( F^d_{n}(t)\) converges to \(F(t) \) for every $ t \in \mathbb{R}^{+} $. Furthermore, the convergence is uniform over \( \mathbb{R}^{+} \) by applying similar arguments as in the proof of Lemma 2.11 in \cite{van2000asymptotic}, using that $F$ is monotone (by non-negativity of $f$), continuous and bounded. This completes the proof of the first claim.

		Now,  $F_n^d$ is continuous on $\mathbb R^+$ so $\lim_{t\to\infty}F_n^d(t)$ exists in $[-\infty,\infty]$ and we have
		\begin{eqnarray*}
			|\lim_{t\to\infty}F_n^d(t)-\lim_{t\to\infty}F(t)|
			&\leq& \Vert F_n^d-F\Vert_\infty=o_p(1),
		\end{eqnarray*}
		where $\lim_{t\to\infty}F(t)=1$. Therefore, 
		\begin{eqnarray}\label{eq: lim F_n^d}
			\lim_{t\to\infty}F_n^d(t)
			=1+o_p(1),
		\end{eqnarray}
		where the right-hand side is finite and positive with probability that tends to one as $n\to\infty$. This proves the second claim of the lemma.
		\hfill{$\Box$}
	}

	\subsection{Proof of Theorem \ref{theo: cvunif}}
	
	{
		We prove the first claim for $f_{n0}$ only, since the proof for $f_0$ is similar. If there exists $x>0$ where $f_{n0}$ is not continuous, then the left-hand slope of $\mathcal MF_n$ at $x$ is strictly larger than its right-hand slope. This implies that $\mathcal MF_n(x)=F_n(x)$ and because $\mathcal MF_n$ is above $F_n$ by definition, this implies that the left-hand slope of $F_n$ at $x$ is strictly larger than its right-hand slope, which contradicts differentiability of $F_n$. Hence, $f_{n0}$ is continuous.
		
		We now prove that $f_{n0}$ is non-negative. Equivalently we prove that $\mathcal M F_n$ is non-decreasing. Arguing by contradiction, we assume that $\mathcal MF_n$ is not non-decreasing. Then there exists $y_n>0$ where the slope $-u_n$ of $\mathcal MF_n$  is stricly negative,  whence $u_n>0$. By concavity of $\mathcal M F_n$ we then have
		\begin{eqnarray*}
			\mathcal M F_n(y_n+t)\leq  \mathcal M F_n(y_n)-u_nt
		\end{eqnarray*}
		for all $t>0$ whence $\lim_{t\to\infty} \mathcal M F_n(t)=-\infty$. Since $\mathcal MF_n$ is above $F_n$ by definition, this implies that $\lim_{t\to\infty} F_n(t)=-\infty$ which contradicts the assumption that $F_n$ is bounded, whence $f_{n0}$ is non-negative.

		We turn to the proof of the second claim. For all $t\geq 0$,  one has 
		$$
		|F_{n}(t) -F(t) | \leq \int_{0}^{t}|f_n(u) -f(u) |du \leq t\sup_{u\in\mathbb{R}^{+}}|f_{n}(u) -f(u) |,
		$$
		where by assumption, the supremum on the right-hand side tends to zero.
		Then, for every $ t \in \mathbb{R}^{+} $,  \( F_{n}(t)\) converges to \(F(t) \). Furthermore, the convergence is uniform over \( \mathbb{R}^{+} \) by applying similar arguments as in the proof of Lemma 2.11 in \cite{van2000asymptotic}, using that $F$ is monotone (by non-negativity of $f$), continuous and bounded. Then, it follows from Marshall's inequality \cite[Exercise 3.1]{groeneboom2014} that 
		\begin{equation*}
			\| \mathcal MF_n- \mathcal MF \|_\infty \leq \| F_{n}-F \|_\infty,
		\end{equation*}
		whence the left-hand side converges to zero.

		Now, we  show that  $f_{n0}$ uniformly converges to $f_0$. It follows from the max-min formulas that for arbitrary $x> 0$ one has
		\begin{eqnarray*}
			f_{n0}(x)=\inf_{v>x}\sup_{u<x}\left(\frac 1{v-u}\int_u^vf_n(t)dt\right)\mbox{ and } f_0(x)=\inf_{v>x}\sup_{u<x}\left(\frac 1{v-u}\int_u^vf(t)dt\right),
		\end{eqnarray*}
		and for all $u<v$. Morever, the triangle inequality yields that
		\begin{eqnarray*}
			\left|\frac 1{v-u}\int_u^vf_n(t)dt-\frac 1{v-u}\int_u^vf(t)dt\right|
			&\leq&\frac 1{v-u}\int_u^v|f_n(t)-f(t)|dt\\
			&\leq& \|f_n-f\|_\infty=o(1).
		\end{eqnarray*}
		Hence,
		\begin{eqnarray*}
			f_{n0}(x)&=&\inf_{v>x}\sup_{u<x}\left(\frac 1{v-u}\int_u^vf_n(t)dt\right)\\
			&=&\inf_{v>x}\sup_{u<x}\left(\frac 1{v-u}\int_u^vf(t)dt\right)+o(1)\\
			&=&f_0(x)+o(1).
		\end{eqnarray*}
		For $x=0$ the min max formulas become
		\begin{eqnarray*}
			f_{n0}(0)=\inf_{v>0}\left(\frac 1{v}\int_0^vf_n(t)dt\right)\mbox{ and } f_0(0)=\inf_{v>0}\left(\frac 1{v}\int_0^vf(t)dt\right),
		\end{eqnarray*}
		and similar as above, the triangle inequality implies that for all $v>0$ one has
		\begin{eqnarray*}
			\left|\frac 1{v}\int_0^vf_n(t)dt-\frac 1{v}\int_0^vf(t)dt\right|
			&\leq& \|f_n-f\|_\infty=o(1).
		\end{eqnarray*}
		Hence,
		\begin{eqnarray*}
			f_{n0}(0)&=&\inf_{v>0}\left(\frac 1{v}\int_0^vf_n(t)dt\right)\\
			&=&f_0(0)+o(1).
		\end{eqnarray*}
		
		We have proved that $f_{n0}(x)$ converges to $f_0(x)$ for all $x\geq 0$. The uniform convergence of $f_{n0}$ follows by continuity, boundedness and monotonicity of $f_0$ with similar arguments as in the proof of Lemma 2.11 in \cite{van2000asymptotic}.

		Consider now the stochastic setting. If   
		$\| f_n-f\|_\infty$  converges to zero in probability, then  from every subsequence one can extract a further subsequence along which $\| f_n-f\|_\infty$ converges to zero almost surely. Along such a subsequence, we know that  $\| f_{n0}-f\|_\infty$ and $\| \mathcal MF_n-\mathcal MF\|_\infty$ converge to zero almost surely. Hence from every subsequence one can extract a further subsequence along which $\| f_{n0}-f_0\|_\infty$ and $\| \mathcal M F_n-\mathcal MF\|_\infty$ converge  to zero almost surely, whence in probability; whence the convergences hold along the whole sequence in probability. This completes the proof of Theorem \ref{theo: cvunif}.
		\hfill{$\Box$}
	}
	
	\subsection{Proof of Corollary \ref{cor}}
	{
		It follows from Lemma \ref{lem: welldef} that $
		\|f_n^d-f\|_\infty=o_p(1)
		$ and $F_n^d\in\ell_\infty(\mathbb R^+)$ with probability that tends to one.
		Hence, Theorem \ref{theo: cvunif} ensures the announced results when  $\mathcal M_n$ is $\mathcal M F_n^d$. The results when $\mathcal M_n$ is $\mathcal M \hat F_n^d$  easily follow using \eqref{eq: lim F_n^d} and \eqref{eq: genuineDF2}. \hfill{$\Box$}

		\subsection{Proof of Theorem \ref{theo: limitdistrib}}\label{sec: prooflimitdistrib}
		
		We have $F_n^d=\tilde F_n^d-\tilde F_n^d(0)$ by definition, and $F(0)=0$, so it follows from  Theorem \ref{theo: ST} that
		\begin{eqnarray*}
			\sqrt n (F_n^d-F) &=&\sqrt n (\tilde F_n^d-F)  -\sqrt n (\tilde F_n^d-F)(0)\\
			&\rightsquigarrow & \mathbb{G}_F-\mathbb{G}_F(0)\quad in\ \ell_\infty(\mathbb R).
		\end{eqnarray*}
		Combining with  the functional delta method and Theorem \ref{theo: brendan} yields
		\begin{equation*}
			\sqrt{n}( \mathcal{M} F_{n}^{d}-\mathcal{M}F) \rightsquigarrow 
			\mathcal{M}_{F}^{\prime }(\mathbb{G}_F-\mathbb{G}_F(0))\ \ in\ \ell_\infty(\mathbb{R}),
		\end{equation*}%
		as $n\to\infty$. The result follows  since for arbitrary function $G\in C_0$ and constant $c\in\mathbb R$ one has  $\mathcal{M}_{F}^{\prime }(G+c)=\mathcal{M}_{F}^{\prime }G+c$. \hfill{$\Box$}
	}

	\section{Proof of Theorem \ref{theo: limit}}\label{sec: prooflimit}
	Define the operator $\mathcal{W}=\mathcal{M}-\mathcal{I}$ where $\mathcal{I}$
	stands for the identity operator. 
	Assume $F$ is concave on $\mathbb{R}^+$. Then $\mathcal{M}F=F$, so
	the statistic in (\ref{Tn}) becomes 
	\begin{equation}\label{eq: Tn}
		T_n(F_{n}^{d})=\sqrt{n}\Vert \mathcal{W}F_{n}^{d} - 
		\mathcal{W}F \Vert_\infty.\end{equation}
	Now, for an arbitrary function $\theta\in\ell_\infty(\mathbb R^+)$ and $a\in\mathbb R$ one has  $\mathcal M(\theta+a)=\mathcal M \theta+a$, which implies that $\mathcal W(\theta+a)=\mathcal W\theta$. Since by definition one has $F^d_n=\tilde F_n^d-\tilde F_n^d(0)$ where $\tilde F_n^d$ is taken from \eqref{tildeFn}, one gets
	\begin{equation}\label{eq: Tntilde}
		T_n(F_{n}^{d})=\sqrt{n}\Vert \mathcal{W}\tilde F_{n}^{d}  -
		\mathcal WF \Vert_\infty.
	\end{equation}

	Now, from Theorem \ref{theo: brendan}, we can deduce that $\mathcal{W}$ is also Hadamard directionally
	differentiable at $F$ (given that $F$ is assumed to be concave) and its derivative is $\mathcal{W}_{F}^{\prime }=\mathcal{M}_{F}^{\prime }-\mathcal{I}$ since the Hadamard directional derivative of the identity is itself. Hence, combining Theorem \ref{theo: ST} with  the functional delta method {and Theorem \ref{theo: brendan}} yields
	\begin{equation*}
		\sqrt{n}( \mathcal{W}\tilde F_{n}^{d}-\mathcal{W}F) \rightsquigarrow 
		\mathcal{M}_{F}^{\prime }\mathbb{G}_F-\mathbb{G}_F\ \ in\ \ell_\infty(\mathbb{R}),
	\end{equation*}%
	as $n\to\infty$.
	Since the sup-norm is a continuous function on $\ell_\infty(\mathbb{R})$, 
	invoking the continuous mapping theorem completes the proof of  the first claim in Theorem \ref{theo: limit}. 
	
	Now, assume that $F$ is not concave on $\mathbb{R}^+$ and define 
	\begin{equation*}
		A:= \Vert F-\mathcal MF\Vert_\infty,
	\end{equation*}
	whence $A$ is strictly positive. It follows from the triangle inequality that
	\begin{eqnarray*}
		n^{-1/2}T_n(F_{n}^{d})&\geq& \Vert F-\mathcal MF\Vert_\infty- \Vert F-F_n^d\Vert_\infty -\Vert \mathcal MF_n^d-\mathcal MF\Vert_\infty.
	\end{eqnarray*}
	Using Marshall's inequality \cite[Exercise 3.1]{groeneboom2014} on gets
	\begin{eqnarray*}
		n^{-1/2}T_n(F_{n}^{d})&\geq& \Vert F-\mathcal MF\Vert_\infty-2 \Vert F-F_n^d\Vert_\infty\\
		&\geq& A- 2\Vert F-F_n^d\Vert_\infty.
	\end{eqnarray*}
	Therefore,
	\begin{eqnarray*}
		\mathbb{P}\left(n^{-1/2}T_n(F_{n}^{d})\leq A/2\right)&\leq& \mathbb{P}\left(\Vert F-F_n^d\Vert_\infty>A/4\right).
	\end{eqnarray*}
	Now, it follows from {Lemma \ref{lem: welldef} that} $\Vert F-F_n^d\Vert_\infty$ converges in probability to zero  and therefore,
	\begin{equation*}
		\lim_{n\to\infty}\mathbb{P}\left(\Vert F-F_n^d\Vert_\infty>A/4\right)=0.
	\end{equation*}
	The second claim in Theorem \ref{theo: limit} follows from the previous two displays by choosing $C=A/2$. 
	\hfill{$\Box$}

	\section{Proofs for Section  \ref{sec: bootstrap}}\label{sec: proofbootstrap}

	\subsection{Proof of Lemma \ref{lem: genuineDF}}
	
	{It follows from the second claim of Lemma \ref{lem: welldef} that $\hat F_n^d$ is well defined with $\hat F_n^d(t)=0$ for all $t\leq 0$ and
		\begin{eqnarray*}
			\lim_{t\to\infty}\hat F_n^d(t)=1
		\end{eqnarray*}
		by definition. Moreover, $\hat F_n^d$ is continuous by concavity and non-decrea-sing by Theorem \ref{theo: cvunif}, whence the lemma follows. \hfill{$\Box$}
	}
	
	\subsection{Proof of Theorem \ref{theo: bootstrap}}
	Recall that from Example 1 in  \cite{sohl2012uniform}, we know that we can choose $\zeta=\mathds{1}_{(-\infty, 0]}$ in Theorem \ref{theo: donsker} provided that $\gamma_s<1/2$ and in that case, $\vartheta_t^*$ is the common distribution function of the iid observations while $\widehat\vartheta_t^*$ is the corresponding deconvolution estimator based on the bootstrap sample.
	Recall moreover that thanks to Corollary \ref{cor} we know that $\hat f_n$ uniformly converges to $f_0$.
	In order to prove Theorem  \ref{theo: bootstrap}, we will need to apply Theorem \ref{theo: donsker} with the above choice of $\zeta$, and {$f$ replaced by $f_0$ in Assumption {\bf B}}, and also  Lemma \ref{lem: genuineDF} so the first step is to show that the assumptions of theorems \ref{theo: ST}  and \ref{theo: donsker} are satisfied under the assumptions of Theorem  \ref{theo: bootstrap}. In the following Lemma, which proof is given in Section \ref{sec: DE}, we show that assumptions of Theorem \ref{theo: ST} hold.

	\begin{lemma}\label{lem: DE} 
		\begin{enumerate}
			\item
			Assumptions {\bf D'} and {\bf E'} imply Assumptions {\bf D} and {\bf E}.
			\item Let $\gamma$ and $\beta$ be positive numbers such that $\gamma>\beta$, and $h_n$ a sequence of positive numbers that converges to zero as $n\to\infty$. If $h_n^{2+2\beta}n\to\infty$ then $h_n^\rho n\to\infty$ for some $\rho>4\beta-4\gamma+2$.
			\item The kernel \( K \) given by \eqref{eq:kernel} 
			satisfies conditions {\bf K} for all $L<rs$. 
		\end{enumerate}
	\end{lemma}

	In the sequel we denote by $\hat f_n$ the slope of the distribution function 
	$\mathcal M \hat F_n^d$ defined by \eqref{eq: genuineDF2} and by $f_n$ the slope of $\mathcal M F_n^d$, which are well defined almost everywhere by concavity. Note that thanks to Lemma \ref{lem: genuineDF}, $\hat f_n$ is a genuine density function (with probability that tends to one)  in the sense that it is non-negative and integrates to one, which is not necessarily the case of $ f_n$. Since we are concerned with convergence in probablity, we can restrict ourselves to a sequence of  events whose probability tends to one as $n\to\infty$ so we assume in the sequel without loss of generality that $\hat f_n$ is a genuine density function. Hence, it follows from Theorem  \ref{theo: donsker} (where we set $\zeta=\mathds{1}_{(-\infty,0]}$ so that  $\vartheta_t^*=\mathcal M \hat F_n^d(t)$ and $\hat\vartheta_t^*=\widetilde F_m^{d*}(t)$, see also Example 1 in \cite{sohl2012uniform}) that the first claim of  Theorem \ref{theo: bootstrap}  holds provided that Assumption {\bf B} is satisfied with $f$ replaced by $f_0$. This means that the first claim is a direct consequence of Corollary \ref{cor}, Lemma \ref{lem: DE}  and the following lemma, that is proved in Section \ref{sec: proofhypfn}.

	\begin{lemma}\label{lem: hypfn}
		Under the assumptions of Theorem \ref{theo: bootstrap}, with $\hat f_n$ the slope of the distribution function 
		$\mathcal M \hat F_n^d$ defined by \eqref{eq: genuineDF2}, one has \eqref{eq: fnhatbounded}, \eqref{eq: fnhatmoments} and \eqref{eq: assumption fii}. 
	\end{lemma}

	We turn to the proof of the second claim  in Theorem \ref{theo: bootstrap}, hence we assume that $F$ is concave. As in Section \ref{sec: prooflimit}
	we consider the operator $\mathcal{W}=\mathcal{M}-\mathcal{I}$ where $\mathcal{I}$
	stands for the identity operator. 
	Since $\hat F_n^d$ is concave on $\mathbb{R}^+$, we have  $\mathcal{M}\hat F_n^d=\hat F_n^d$, so the bootstrap version $T_m(F_m^{d*})$ of
	the test statistic  becomes 
	\[
	T_m(F_m^{d*})=\sqrt{m}\Vert \mathcal{W}F_{m}^{d*} - 
	\mathcal{W}\hat F_n^d \Vert_\infty.\] 
	
	{It follows from  Theorem \ref{theo: ST} combined with the continuous mapping theorem that
		\begin{equation*}
			\sqrt n\sup_{t\geq 0}|F_n^d(t)-F(t)|\stackrel{L}{\longrightarrow}\sup_{t\geq 0}|\mathbb{G}_F(t)-\mathbb{G}_F(0)|
		\end{equation*}
		as $n\to\infty$ with $\mathbb{G}_F$ a centered Gaussian Borel random variable in $\ell_\infty(\mathbb{R})$, whence the left-hand side is bounded in probability. Hence,  with  $(\epsilon_n)_n$ an arbitrary sequence of positive numbers that converges to zero as $n\to\infty$ with $\lim_{n\to\infty}\sqrt n\epsilon_n=+\infty$, one has
		\begin{eqnarray}\label{eq: supeps}
			\Vert F_n^d-F\Vert_\infty\leq \epsilon_n
		\end{eqnarray}
		with probability that tends to one. Moreover,
		\begin{eqnarray*}
			\left\vert\lim_{t\to\infty}F_n^d(t)-1\right\vert
			&=&\left\vert\lim_{t\to\infty}F_n^d(t)-
			\lim_{t\to\infty}F(t)\right\vert\\
			&\leq&\|F_n^d-F\|_\infty
		\end{eqnarray*}
		whence we have
		\begin{eqnarray}\label{eq: denom}
			\left\vert\lim_{t\to\infty}F_n^d(t)-1\right\vert
			&\leq&\epsilon_n
		\end{eqnarray}
		with probability that tends to one. Combining} with the fact that the test statistic in \eqref{eq: Tn} is bounded in probability where $\mathcal WF=0$, one gets
	\[
	{\sqrt m \Vert \mathcal{W}F_n^d  - 
		\mathcal{W}\hat F_n^d \Vert_\infty=O_{\mathbb P}(m^{1/2}n^{-1/2}\epsilon_n).}\] 
	Since $m\ll n$ one then gets
	\begin{eqnarray*}
		T_m(F_m^{d*})=\sqrt{m}\Vert \mathcal{W}F_{m}^{d*} - 
		\mathcal{W} F_n^d \Vert_\infty+o_{\mathbb P}(1)\\
		=\sqrt{m}\Vert \mathcal{W}F_{m}^{d*} - 
		\mathcal{W} F \Vert_\infty+o_{\mathbb P}(1),
	\end{eqnarray*}
	using again that  the test statistic in \eqref{eq: Tn} is bounded in probability for the last equality.  Since for arbitrary functions $\theta\in\ell_\infty(\mathbb R^+)$ and $a\in\mathbb R$ one has  $\mathcal M(\theta+a)=\mathcal M \theta+a$, which implies that $\mathcal W(\theta+a)=\mathcal W\theta$, we also have
	\[
	T_m(F_m^{d*})=\sqrt{m}\Vert \mathcal{W}\widetilde F_{m}^{d*} - \mathcal W F \Vert_\infty+o_{\mathbb P}(1).\] 
	In fact, the assumption that that $m\ll n$ is made precisely with the aim to replacing $F_n^d$ by $F$ in the above displays, and also in \eqref{eq: cvboot} which, thanks to \eqref{eq: ST} and Marshall's inequality,  can be reformulated as
	\begin{equation*}
		\sqrt m( \widetilde F_m^{d*}-F) \rightsquigarrow \mathbb{G}_F\ \ in\ \ell_\infty(\mathbb{R}),
	\end{equation*}%
	as $n\to\infty$. 
	Now, we have seen in Section \ref{sec: prooflimit} that  $\mathcal{W}$ is Hadamard directionally
	differentiable at $F$ with  derivative $\mathcal{W}_{F}^{\prime }=\mathcal{M}_{F}^{\prime }-\mathcal{I}$ so one gets
	\begin{equation*}
		\sqrt{m}( \mathcal{W}\widetilde F_{m}^{d*}-\mathcal{W} F) \rightsquigarrow 
		\mathcal{M}_{F}^{\prime }\mathbb{G}_F-\mathbb{G}_F\ \ in\ \ell_\infty(\mathbb{R}),
	\end{equation*}%
	as $n\to\infty$, conditionally. 
	Since the sup-norm is a continuous function on $\ell_\infty(\mathbb{R})$, 
	invoking the continuous mapping theorem completes the proof of  the 
	second claim  in Theorem \ref{theo: bootstrap}.
	
	It remains to prove the third claim. It follows from the triangle inequality that
	\begin{eqnarray*}
		\Vert \mathcal{M}F_{m}^{d*} -F_{m}^{d*} \Vert_\infty
		&\leq & \Vert \mathcal{M}F_{m}^{d*} - \mathcal{M} \hat F_{n}^{d} \Vert_\infty
		+\Vert \mathcal{M} \hat F_{n}^{d}  -F_{m}^{d*} \Vert_\infty
	\end{eqnarray*}
	where $\mathcal{M} \hat F_{n}^{d} =\mathcal M\mathcal{M} \hat F_{n}^{d} $ by concavity. 
	Marshall's inequality \cite[Exercise 3.1]{groeneboom2014} then  yields
	\begin{eqnarray*}
		\Vert \mathcal{M}F_{m}^{d*} -F_{m}^{d*} \Vert_\infty
		&\leq&\Vert F_{m}^{d*} - \mathcal{M} \hat F_{n}^{d} \Vert_\infty
		+\Vert \mathcal{M} \hat F_{n}^{d}  - F_{m}^{d*} \Vert_\infty\\
		&=  & 2\Vert  F_{m}^{d*}-\mathcal{M} \hat F_{n}^{d}  \Vert_\infty\\
		&\leq&2\Vert  \widetilde F_{m}^{d*}-\mathcal{M} \hat F_{n}^{d}  \Vert_\infty+2\vert  \widetilde F_{m}^{d*}(0) \vert
	\end{eqnarray*}
	Now, it follows from the continuous mapping theorem combined to Theorem \ref{theo: bootstrap} that $\sqrt m\Vert \widetilde F_{m}^{d*}-\mathcal{M} \hat F_{n}^{d}  \Vert_\infty$ converges in distribution to $\Vert \mathbb G_{F_0} \Vert_\infty$ while $\sqrt m \vert  \widetilde F_{m}^{d*}(0) \vert$ converges in distribution to $| \mathbb G_{F_0}(0)|$, conditionally as $n\to\infty$. This implies that conditionally, the variable in the last line of the previous display  is bounded in probability. This implies that $\sqrt  m\Vert \mathcal{M}F_{m}^{d*} -F_{m}^{d*} \Vert_\infty$ is dominated by a variable that is bounded in probability, hence it is bounded in probability as well, conditionally. This completes the proof of Theorem \ref{theo: bootstrap}.
	\hfill{$\Box$}

	\subsection{Proof of Theorem \ref{lem: level}}\label{sec: prooflevel}

	In the sequel, for arbitrary concave function $F$ on $\mathbb R^+$ we denote by $H_F$ the distribution function of the limit distribution $\|\mathcal M_{F}'\mathbb{G}_{F}-\mathbb{G}_{F}\|_\infty$, and by $r_F$  the left-hand point of its support, defined as the infimum of those $r\geq 0$ for which $H_F(r)>0$.

	It follows from Theorem \ref{theo: ST} and Theorem \ref{theo: bootstrap} that the conditional limit distribution of $T_m(F_m^{d*})$ is the same as the limit distribution of $T_n(F_n^{d})$, so the test defined by the critical region in \eqref{eq: critical} satisfies the significance level constraint \eqref{eq: level} 
	provided that the limiting distribution function $H_F$  is continuous at its $(1-\gamma)$-quantile, $c_{\gamma,F}$ say. (This can be proven with similar arguments as for the proof of Lemma 23.2  in \cite{van2000asymptotic}.) As pointed out in \cite[page 4605]{fang2019refinements} the map $G\mapsto \|\mathcal M'_F G-G\|_\infty$ is convex on $\ell_\infty(\mathbb R^+)$ so Theorem 11.1 in \cite{davydov1998local} applies. It follows from that theorem that $H_F$ is everywhere continuous except possibly at $r_F$. This implies that \eqref{eq: level} holds provided that $c_{\gamma,F}>r_F$ or 
	equivalently, provided that 
	\begin{equation*}
		H_F(r_F)<1-\gamma.
	\end{equation*}
	Since it is assumed that $\gamma<1/2$,  it suffices to prove
	\begin{equation}\label{eq: davydov}
		H_F(r_F)\leq\frac12.
	\end{equation}

	To prove \eqref{eq: davydov}, we first show that $r_F=0$. 
	For arbitrary $r>0$ one has
	\begin{eqnarray*}
		H_F(r)&=&\mathbb P\left(\|\mathcal M'_F\mathbb G_F-\mathbb G_F\|_\infty\leq r\right)
	\end{eqnarray*}
	where by the triangle inequality,
	\begin{eqnarray*}
		\|\mathcal M'_F\mathbb G_F-\mathbb G_F\|_\infty\leq\|\mathcal M'_F\mathbb G_F\|_\infty+\|\mathbb G_F\|_\infty.
	\end{eqnarray*}
	It follows from Theorem \ref{theo: brendan} that for all $t>0$ either $\mathcal M'_F\mathbb G_F(t)=\mathbb G_F(t)$, or $\mathcal M'_F\mathbb G_F$ is affine on an interval that contains $t$ in its interior and such that $\mathcal M'_F\mathbb G_F$ coincides with $\mathbb G_F$ at the boundaries of that interval. This implies that $\|\mathcal M'_F\mathbb G_F\|_\infty\leq\|\mathbb G_F\|_\infty$ whence
	\begin{eqnarray*}
		\|\mathcal M'_F\mathbb G_F-\mathbb G_F\|_\infty\leq2\|\mathbb G_F\|_\infty,
	\end{eqnarray*}
	and therefore,
	\begin{eqnarray}\label{eq: H_F}
		H_F(r)&\geq&\mathbb P\left(\|\mathbb G_F\|_\infty\leq\frac r2\right).
	\end{eqnarray}

	By Theorem~1 and Theorem~7 of \cite{sohl2012uniform}, the limit process
	$\mathbb G_F$ is a centered Gaussian Borel random element in
	$\ell_\infty(\mathbb R^+)$ and admits a separable version with uniformly
	$d$-continuous sample paths and $\sup_{t\ge0}|\mathbb G_F(t)|<\infty$
	almost surely, where
	\[
	d(s,t):=\left(\mathbb E\left[(\mathbb G_F(s)-\mathbb G_F(t))^2\right]\right)^{1/2}
	\]
	is the intrinsic covariance metric. Uniform $d$-continuity
	means that, for almost every $\omega$ and for every $\varepsilon>0$, there exists
	$\delta>0$ such that $d(s,t)<\delta$ implies
	$|\mathbb G_F(s,\omega)-\mathbb G_F(t,\omega)|<\varepsilon$.
	Moreover, Theorem~7 of \cite{sohl2012uniform} establishes that
	$\mathbb R^+$ is totally bounded with respect to $d$, meaning that, for every $\varepsilon>0,$ $\mathbb R^+$ can be covered by finitely many balls of radius $\varepsilon$ \citep[p. 17]{aadbook}, and that the above version of
	$\mathbb G_F$ is a tight Borel random element in $\ell_\infty(\mathbb 	R^+)$. 
	According to Section 3.1 in \cite{bogachev1998} (page 97), a measure on a complete, separable, metrizable and localy convex space is a Radon measure.
	Now, define the space $B$ of bounded and uniformly $d$-continuous real-valued functions on $\mathbb R^+$,
	endowed with the supremum norm. The set of uniformly continuous functions on a totally bounded, semimetric
	space is complete and separable \citep[p. 38]{aadbook}. Since $\mathbb R^+$
	is totally bounded with respect to $d$, then $B$ is complete and separable, and the separable version of
	$\mathbb G_F$ takes values in $B$ almost surely. Hence the law
	of $\mathbb G_F$, denoted as $\gamma,$ is a centered Gaussian probability measure on the
	separable, locally convex space $B$; tightness and Borel measurability
	imply that $\gamma$ is a Radon measure on $B$.
	By Theorem~3.6.1 in \cite{bogachev1998}, the topological support of a
	Gaussian Radon measure on a locally convex space coincides with
	$a_\gamma+\overline{H(\gamma)}$, where $a_\gamma$ denotes the mean and $H(\gamma)$ the
	Cameron-Martin space. Note that we don't need to characterize $\overline{H(\gamma)}$ in our case, it suffices to note that $H(\gamma)$ \citep[p. 44]{bogachev1998} contains 0 (the function which is identically 0).
	Since $\gamma$ is centered we have $a_\gamma=0$, hence
	$0$ belongs to the support of $\gamma$. By the definition of topological support,
	every open neighborhood of $0$ in $B$ has strictly positive $\gamma$-measure.
	In particular, for every $r>0$,	$\gamma\bigl(\{g\in B:\|g\|_\infty<r\}\bigr)>0,$ that is, $$
	P(\|\mathbb G_F\|_\infty<r)>0.$$
	Combining this with \eqref{eq: H_F} yields $H_F(r)>0$ for all $r>0$, and thus $r_F=0$.

	We turn to the proof of \eqref{eq: davydov}. For this task, we denote by $I$ an interval on which $F$ is affine. 
	Since $r_F=0$  we have
	\begin{equation*}
		H_F(r_F)=\mathbb P\left(\|\mathcal M'_F\mathbb G_F-\mathbb G_F\|_\infty=0\right).
	\end{equation*}
	If $\|\mathcal M'_F\mathbb G_F-\mathbb G_F\|_\infty$ is equal to zero, then $\mathcal M'_F\mathbb G_F=\mathbb G_F$ on $\mathbb R^+$. In particular, the equality holds on $I$ and because $\mathcal M'_F\mathbb G_F$ is concave on $I$, this implies that $\mathbb G_F$ is also concave on $I$. Hence,
	\begin{eqnarray*}
		H_F(r_F)
		&\leq&\mathbb P\left(\mathbb G_F\mbox{ is concave on }I\right).
	\end{eqnarray*}
	Since $F$ is affine on  $I$, the process $\mathbb G_F$ is concave on $I$ if and only if $\mathbb H_F:=(\mathbb G_F(t)-ZF(t))_{t\in\mathbb R}$ where $Z$ is an independant standard Gaussian variable, is concave on $I$, whence
	\begin{eqnarray*}
		H_F(r_F)
		&\leq&\mathbb P\left(\mathbb H_F\mbox{ is concave on }I\right).
	\end{eqnarray*}
	Note that  $\mathbb H_F$ is a centered Gaussian process with covariance function 
	$\Sigma'_{s,t}=\int g_tg_sf_Y$, which is slightly simpler than that of $\mathbb G_F$. 
	
	Fix any two points $x<z$ in $I$ and for every $y\in(x,z)$ set
	\[
	\alpha_y := \frac{z-y}{z-x}\in(0,1); \qquad 
	L_y :=\mathbb{H}_F(y)- \alpha_y\mathbb{H}_F(x)-(1-\alpha_y)\mathbb{H}_F(z)
	\]
	and $V(y)=\operatorname{Var}(L_y)$. 
	We show below that  $V$ is not identically zero on $(x,z)$. Otherwise said, with probability one $\mathbb H_F$ (as well as $\mathbb G_F$) cannot be affine on the interval.
	Hence, we can find a particular $y_0\in I$ such that $V(y_0)>0$,  whence the random variable $L_{y_0}$ is non-degenerate.
	If $\mathbb H_F$ is concave on $I$ then $L_y\geq 0$ for all $y\in(x,y)$, and in particuler for $y=y_0$.  Since $L_{y_0}$ is non-degenerate centered Gaussian, we arrive at
	\[
	H_F(r_F)\leq\mathbb P\left(\mathbb H_F\mbox{ is concave on }I\right)\leq \mathbb{P}(L_{y_0}\ge 0)=\frac12
	\]
	which completes the proof of \eqref{eq: davydov} and hence, completes the proof of Lemma \ref{lem: level}. 
	
	It remains to show that $V$ is not identically zero on $(x,z)$. We have 
	\begin{eqnarray*}
		L_y  &=&\alpha_y\left(\mathbb{H}_F(y)-\mathbb{H}_F(x)
		\right)+(1-\alpha_y)\left(\mathbb{H}_F(y)-\mathbb{H}_F(z)
		\right)
	\end{eqnarray*}
	and therefore,
	\begin{eqnarray*}
		V(y)&=&\alpha_y^2\int f_Y(g_y-g_x)^2+2\alpha_y(1-\alpha_y)\int f_Y(g_y-g_z)(g_y-g_x)\\&&+(1-\alpha_y)^2\int f_Y(g_y-g_z)^2\\
		&=&\int f_Y\left(\alpha_y(g_y-g_x)+(1-\alpha_y)(g_y-g_z)\right)^2.
	\end{eqnarray*}
	Arguing by contradiction, let us assume that $V(y)=0$ for all $y\in(x,z)$. The above display implies that on the support of $f_Y$,
	\begin{eqnarray*}
		\alpha_y(g_y-g_x)+(1-\alpha_y)(g_y-g_z)\equiv 0
	\end{eqnarray*}
	where $g_u=h\star \mathds{1}_{(-\infty, u]}$ for all $u$, and
	\begin{equation*}
		h= \mathcal{F}^{-1} \left[ \frac{1}{\mathcal{F}f_\varepsilon(-\bullet)} \right]. 
	\end{equation*}
	Hence we have
	\begin{eqnarray*}
		\alpha_y\int_x^y h(t-v)dv-(1-\alpha_y)\int_y^zh(t-v)dv= 0
	\end{eqnarray*}
	for all $t$ in the support of $f_Y$, and all $y\in (x,z)$.

	The function being constant on an interval, its derivative with respect ot $y$ is constant equal to zero almost everywhere on that interval whence 
	\begin{eqnarray*}
		-\frac{1}{z-x}\int_x^yh(t-v)dv-\alpha_yh(t-y)-\frac{1}{z-x}\int_y^zh(t-v)dv-(1-\alpha_y)h(t-y)=0
	\end{eqnarray*}
	for all $t$ in the support of $f_Y$ and almost all $y\in (x,z)$. We arrive at 
	\begin{eqnarray*}
		-\frac{1}{z-x}\int_x^zh(t-v)dv=h(t-y)
	\end{eqnarray*}
	for all $t$ in the support of $f_Y$ and almost all $y\in (x,z)$. The left-hand side does not depend on $y$ so for a fixed $t$, the function $y\mapsto h(t-y)$ is almost everywhere constant on $(x,z)$. 
	If we assume that the support of $\varepsilon$ is $\mathbb R$, then the support of $Y$ is $\mathbb R$ as well, and the above property implies that $h$ is constant almost everywhere on the whole real line. This means that $h$ is the inverse Fourier of the function $c\delta_0$ for some constant $c\in\mathbb R$, where for all $t\in\mathbb R$, $\delta_0(t)$ is equal to one if $t=0$ and to zero otherwise. By unicity of the inverse Fourier transform (in the distributional sense) this implies that 
	$$\frac{1}{\mathcal{F}f_\varepsilon(-\bullet)} =c\delta_0.$$
	This is not possible since the left-hand side cannot equal zero contrary to the right-hand side, whence $V$ cannot be identically zero on $(x,z)$. \hfill{$\Box$}

	\subsection{Proof of Lemma \ref{lem: DE}}\label{sec: DE}
	Suppose that Assumptions {\bf D'} and {\bf E'} hold. Then there exists a positive $C$ such that 	
	$\mathbb{E}(X^{8}) \leq C$ and $\mathbb{E}(\varepsilon^{8}) \leq C$. By Jensen's inequality, this implies that $\int \left\vert x\right\vert ^{2+\delta }f\left( x\right) dx<\infty $ and $\int \left\vert x\right\vert ^{2+\delta }f_{\varepsilon
	}\left( x\right) dx<\infty $ for  all  positive $\delta \leq 6.$ Hence, {\bf E}  and  {\bf D} hold.
	
	The second claim follows from the fact that $2+2\beta>4\beta-4\gamma+2$ for all $\gamma>\beta$.
	
	We turn to the third claim. By construction, \( \mathcal{F}K \) is compactly supported in \( [-1,1] \), which implies that \( K \) is band-limited. Since \( \mathcal{F}K \) is even, \( K \) is real-valued and symmetric. Moreover, as a compactly supported polynomial, \( \mathcal{F}K \in C_c^\infty(\mathbb{R}) \), the space of smooth functions with compact support, and hence its inverse Fourier transform \( K \) belongs to the Schwartz space \( \mathcal{S}(\mathbb{R}) \), defined as the space of all functions that are infinitly differentiable  and such that for all non-negative integers $m$ and $n$, the $n$th derivative $f^{(n)}$ is such that the function $x\mapsto|x^mf^{(n)}(x)|$ is bounded. This implies that $K \in L^{1}(\mathbb{R}) \cap \ell_{\infty }(\mathbb{R})$ whence the first part of assumption {\bf K} is satisfied. Now, we use the identity 
	\[ \int x^j K(x)\, dx = (-i)^j \left. \frac{d^j}{dt^j} \mathcal{F}K(t) \right|_{t=0} ,\] 
	which holds since \( \mathcal{F}K \in C_c^\infty(\mathbb{R}) \).  By construction, the right-hand side is equal to zero  for all \( j = 1, \dots, rs - 1 \), and thus \( \int x^j K(x)\, dx = 0 \) for \( j = 1, \dots, L \), for any integer \( L < rs \). Since \( K\in\mathcal{S}(\mathbb{R}) \), all the absolute moments of $K$ are finite so the second part of assumption {\bf K} holds, as well as the third part, using that $\int K=\mathcal FK(0)=1$. \hfill{$\Box$}\\

	\subsection{Proof of Lemma \ref{lem: hypfn}}\label{sec: proofhypfn}
	In the sequel we denote by $f_n^d$ the slope of the deconvolution estimator $F_n^d$ and by $\hat f_n^d$ the slope of  $\hat F_n^d$ defined by \eqref{eq: genuineDF}, whence
	\begin{eqnarray}\label{eq: ftildetofhat}
		\hat f_n^d(t)=\frac{f_n^d(t)}{\lim_{t\to\infty}F_n^d(t)}
	\end{eqnarray}
	for all $t$. We prove in Sections \ref{sec: boundarydeconv} and \ref{sec: momentdeconv} below that under the assumptions of Theorem \ref{theo: bootstrap} one has
	\begin{equation}\label{eq: boundarydeconv}
		\sup_{t\in\mathbb R^+}|f_n^d(t)|=O_p(1),
	\end{equation}
	and because $r>4$ there exists $\delta>0$ such that
	\begin{equation}\label{eq: momentdeconv}
		\int_0^\infty |t|^{2+\delta}f_n^d(t)dt
		=O_p(1).
	\end{equation}
	Using those properties of the estimator, the equations  \eqref{eq: fnhatbounded}, \eqref{eq: fnhatmoments} and \eqref{eq: assumption fii} are proved in the subsequent sections \ref{sec: fnhatbounded}, \ref{sec: fnhatmoments} and \ref{sec: assumption fii}, respectively. \hfill{$\Box$}

	\subsubsection{Proof of \eqref{eq: boundarydeconv}}\label{sec: boundarydeconv}
	Recall that by assumption,  $f\in H^{\alpha }\left( \mathbb{R}\right) $ for some  $\alpha>1/2$. It follows from the Cauchy-Schwarz's inequality that
	\begin{eqnarray*}
		\int|\mathcal Ff|&=&\int(1+x^2)^{\alpha/2}|\mathcal Ff(x)|\times (1+x^2)^{-\alpha/2}dx\\
		&\leq&\left(\int(1+x^2)^{\alpha}|\mathcal Ff(x)|^2dx\right)^{1/2}\times\left(\int(1+x^2)^{-\alpha}dx\right)^{1/2},
	\end{eqnarray*}
	where the first term on the last line is bounded because  $f\in H^{\alpha }\left( \mathbb{R}\right) $, and the second one if finite since $\alpha>1/2$. Hence, 
	\begin{equation}\label{eq: lapfini}
		\int|\mathcal Ff|<\infty
	\end{equation}
	is finite so we can apply Corollary \ref{cor} to obtain that $\| f_n^d-f\|_\infty=o_p(1)$. Combining with boundedness of $f$ and the triangle inequality proves that $\| f_n^d\|_\infty=O_p(1)$. \hfill{$\Box$}

	\subsubsection{Proof of \eqref{eq: momentdeconv}}\label{sec: momentdeconv}
	Let $\delta>0$ be sufficently small so that $2+\delta<4$. Let $q=2+\delta$, $k=1$,
	\( v_q = \mathbb{E}(|X|^q) \)  and $\hat{v}_q = \int_{-\infty}^{\infty} |t|^q \, f_n^d(t)dt$. Under assumption \eqref{eq: hypphiepsion} one has
	\begin{equation*}
		\left\vert \left( 1/\varphi_{\varepsilon} \right)(t) \right\vert \leq C' (1 + |t|^2)^{\beta /2}, \quad \text{for all } t \in \mathbb{R}.
	\end{equation*}
	Hence, under our Assumptions {\bf D'} and {\bf E'} we have $F_\epsilon\in\mathcal F_5(\beta,C_1)$ and $F_X\in\mathcal G_6(0,C_2)$ in the notations of \cite{hall2008estimation}. Moreover,  the assumption that \( r >4\) implies that \(r> \max(q,4)\)  so it follows from Theorem 3.6 in \cite{hall2008estimation}  that 
	\begin{equation*}
		\mathbb{E}[(\hat{v}_q - v_q)^2] \leq C \left( h^{2q} + n^{-1} h^{-(2\beta - 2q - 1)} \right).
	\end{equation*}
	The upper bound tends to zero since $h\to0 $ and $\beta<1/2$, so it follows from Markov's inequality that $(\hat{v}_q - v_q)^2=o_P(1)$, whence $\hat{v}_q = v_q+o_P(1)$. From the assumption that \( v_{8} \) is finite, we can deduce that \( v_{q}  \) is finite  whence  $\hat{v}_q = O_P(1)$, which proves \eqref{eq: momentdeconv}.
	\hfill{$\Box$}

	\subsubsection{Proof of \eqref{eq: fnhatbounded}}\label{sec: fnhatbounded}
	Since $\hat f_n$ is  the  slope of the least concave majorant of $\hat F_n^d$ and $\hat F_n^d(0)=0$ by definition, we have
	\begin{eqnarray*}
		\hat f_n(0)&=&\sup_{t>0}\frac{\hat F_n^d(t)}{t}.
	\end{eqnarray*}
	Moreover,  it follows from Lemma \ref{lem: genuineDF} that $\hat f_n$ is a genuine density function, in particular, it is non-negative at every point (with probability that tends to one). Since \eqref{eq: denom} holds with probablity that tends to one, it follows from the definition \eqref{eq: genuineDF} that
	\begin{eqnarray}\label{eq: boundfnhat}
		0\leq \hat f_n(0)&\leq &\sup_{t>0}\frac{F_n^d(t)}{t}\times \frac{1}{1-\epsilon_n}
	\end{eqnarray}
	with probability that tends to one. Here,  $(\epsilon_n)_n$ is an arbitrary  sequence of positive numbers that converges to zero as $n\to\infty$ with $\lim_{n\to\infty}\sqrt n\epsilon_n=+\infty$. We will split the above supremum into two parts. We begin with the supremum near zero. Let $T>0$. Since $f_n^d$ is the derivative of $F_n^d$ we have
	\begin{eqnarray*}
		\sup_{t\in(0,T]}\frac{F_n^d(t)}{t}&=&
		\sup_{t\in(0,T]}\left(\frac 1t\int_0^t f_n^d(x)dx\right)\\
		&\leq&\sup_{t\in(0,T]}|f_n^d(t)|=O_p(1),
	\end{eqnarray*}
	where we used \eqref{eq: boundarydeconv} for the last equality.
	We now consider the supremum far from zero. Since \eqref{eq: supeps} holds with probability that tends to one we obtain 
	\begin{eqnarray*}
		\sup_{t\geq T}\frac{F_n^d(t)}{t}
		&\leq&\frac 1T\left(F(t)+\epsilon_n\right)
	\end{eqnarray*}
	with probablity that tends to one, so that the supremum on the left hand side is bounded in probability. The supremum over all $t>0$ is bounded from above by the maximum between the supremum near zero and the supremum far from zero so we obtain that 
	\begin{eqnarray*}
		\sup_{t>0}\frac{F_n^d(t)}{t}=O_p(1).
	\end{eqnarray*}
	Combining with \eqref{eq: boundfnhat} yields that
	$\hat f_n(0)=O_p(1).$
	Now  $\hat f_n$ is a non-increasing function on $[0,\infty)$ as being the slope of a concave function $\mathcal M F_n^d$ whence, 
	\begin{equation*}
		\sup_x\hat f_n(x)\leq \hat f_n(0)=O_p(1).
	\end{equation*}
	This completes the proof of \eqref{eq: fnhatbounded}. \hfill{$\Box$}

	\subsubsection{Proof of \eqref{eq: fnhatmoments}}\label{sec: fnhatmoments}
	The following proof is inspired by Section 1.A.1 in \cite{shaked2007stochastic}. Because the function $g:t\mapsto |t|^{2+\delta}$ is non-decreasing and takes value zero at point zero, it is possible, for each $k\in\mathbb{N}$, to define a sequence of upper sets (that is, of open or  closed right half lines) $U_1,\dots,U_k$, and a sequence of positive numbers $a_1,\dots,a_k$ (all of which may depend on $k$) such that as $k\to\infty$ one has
	\begin{equation*}
		g(t)=\lim_{k\to\infty}\sum_{i=1}^ka_i \mathds{1}_{U_i}(t)
	\end{equation*}
	for all $t$, and the approximating functions $g_k:=\sum_{i=1}^ka_i \mathds{1}_{U_i}$ are such that $g_k\geq g_\ell$ for all $k\geq \ell$. To alleviate the notation we write $M_n=\mathcal M \hat F_n^d$, so that $\hat f_n$ is the slope of $M_n$. Recall that from Lemma \ref{lem: genuineDF}, we can assume without loss of generality that $\hat f_n$ is a probability density function, whence it is non negative and integrates to one. It follows from the monotone convergence theorem that
	\begin{eqnarray*}
		\int |t|^{2+\delta}\hat f_n(t)dt
		&=&\lim_{k\to\infty}\int\sum_{i=1}^ka_i \mathds{1}_{U_i}(t)\hat f_n(t)dt\\
		&=&\lim_{k\to\infty}\sum_{i=1}^ka_i (1-M_n(u_i))
	\end{eqnarray*}
	where $u_i$ denotes the left hand boundary of $U_i$.  But $M_n$ is above $\hat F_n^d$ as being the least concave majorant of the latter function so the previous equalities yield
	\begin{eqnarray*}
		\int |t|^{2+\delta}\hat f_n(t)dt
		&\leq&\lim_{k\to\infty}\sum_{i=1}^ka_i (1-\hat F_n^d(u_i)).
	\end{eqnarray*}
	Now,  the slope of $\hat F_n^d$ is $\hat f_n^d$ that integrates to one, so $(1-\hat F_n^d(u_i))=\int \mathds{1}_{U_i}(t) \hat f_n^d(t)dt$ which yields
	\begin{eqnarray*}
		\int |t|^{2+\delta}\hat f_n(t)dt
		&\leq&\lim_{m\to\infty}\int \sum_{i=1}^ma_i \mathds{1}_{U_i}(t) \hat f_n^d(t)dt\\
		&=&\int g(t) \hat f_n^d(t)dt,
	\end{eqnarray*}
	using again the monotone convergence theorem (possibly decomposing $\hat f_n^d$ into the sum of its positive and negative parts) for the last equality. Combining with \eqref{eq: ftildetofhat} and 
	\eqref{eq: supeps}, that holds with probability that tends to one for some $\epsilon_n\to0$, then yields
	\begin{eqnarray*}
		\int |t|^{2+\delta}\hat f_n(t)dt
		&\leq&\frac 1{1-\epsilon_n}\left|\int g(t)  f_n^d(t)dt\right|
	\end{eqnarray*}
	with probability that tends to one. Combining with \eqref{eq: momentdeconv}  completes the proof of \eqref{eq: fnhatmoments}. \hfill{$\Box$}

	\subsubsection{Proof of \eqref{eq: assumption fii}}\label{sec: assumption fii}
	Let $\alpha\in(0,1/2)$. Since $\hat f_n$ is monotone,  it follows from integration by parts that
	\begin{eqnarray*}
		\mathcal F \hat f_n(x)&=&\left[\frac{e^{iux}}{ix}\hat f_n(u)\right]_0^\infty 
		-\frac{1}{ix}\int_0^\infty e^{iux}d\hat f_n(u).
	\end{eqnarray*}
	Moreover, Lemma \ref{lem: genuineDF} ensures that $\hat f_n$ is non negative and integrates to one with probability that tends to one  so (because one can restrict without loss of generality to an event whose probability tends to one in order to prove boundedness in probability) we can assume without loss of generality that $\lim_{u\to\infty}\hat f_n(u)=0$, whence
	\begin{eqnarray*}
		|\mathcal F \hat f_n(x)|&=&\left|-\frac{1}{ix}\hat f_n(0)
		-\frac{1}{ix}\int_0^\infty e^{iux}d\hat f_n(u)
		\right|\\
		&\leq&\frac{1}{|x|}|\hat f_n(0)|+
		\frac{1}{|x|}\int_0^\infty d\hat f_n(u),
	\end{eqnarray*}
	using that $|e^{iux}|=1$. Hence,
	\begin{eqnarray*}
		|\mathcal F \hat f_n(x)|
		&\leq&\frac{1}{|x|}V_n.
	\end{eqnarray*}
	where $V_n:=|\hat f_n(0)|+1$ can be assumed to be finite thanks to \eqref{eq: fnhatbounded}.
	Since we also have $|\mathcal F \hat f_n(x)|\leq 1$ for all $x$ we conclude that
	\begin{eqnarray*}
		\int\langle u\rangle^{2\alpha}|\mathcal F \hat f_n(u)|^2 du&\leq& 2+\int_{|u|>1} \langle u\rangle^{2\alpha}|\mathcal F \hat f_n(u)|^2 du\\
		&\leq& 2+V_n^2 \int_{|u|>1} \langle u\rangle^{2\alpha}u^{-2} du.
	\end{eqnarray*}
	The integral on the last line is finite since $2\alpha-2<-1$, which completes the proof that \eqref{eq: assumption fii} holds 
	for all $\alpha\in(0,1/2)$.  \hfill{$\Box$}


	\section{Proofs for Section \ref{sec: generalization}}\label{sec: proofdonsker}
	The proof of Theorem \ref{theo: donsker} is given in a first subsection that also states a few lemmas that are required for the proof of the theorem. The lemmas are proved in separate subsections afterwards.
	\subsection{Proof of Theorem \ref{theo: donsker}}
	For notational convenience (to have a single index that tends to infinity), we consider only the case where $m=n$; the general case can be obtained likewise. Moreover, we write $h$ instead of $h_m$ to alleviate the notation.

	Our aim is to show that $
	\pi(\mu_n^*,\mu)$ converges in probability to zero as $n\to\infty$,
	where $\mu_n^*$ is the conditional distribution of $\sqrt n(\widehat \vartheta^*_t- \vartheta_t^*)_{t\in\mathbb R}$, $\mu$ is the distribution of $\mathbb G_F$ and $\pi$ is the Prohorov distance, see the comments following Theorem \ref{theo: bootstrap}. We will prove below that this holds if, instead of assuming that $\|\hat f_n-f\Vert_\infty$ converges to zero in probability, we make the stronger assumption that this convergence holds almost surely. This suffices for the following reason. If $\|\hat f_n-f\Vert_\infty$ converges to zero in probability, then from every subsequence one can extract a further subsequence along which $\|\hat f_n-f\Vert_\infty$ converges to zero almost surely. Along that further subsequence, we know that $
	\pi(\mu_n^*,\mu)$ converges in probability to zero. Hence from every subsequence one can extract a further subsequence along which $
	\pi(\mu_n^*,\mu)$ converges in probability to zero; whence the convergence holds along the whole sequence.

	Hence, in the sequel we assume without loss of generality  that $\|\hat f_n-f\Vert_\infty$ converges to zero almost surely.
	For some fixed $C>0$, we denote by $E_n$ the event
	\begin{eqnarray*}\label{eq: En}\notag
		E_n&:=&\left(\sup_t|\hat f_n(t)|\leq C\right)\cap\left( \int \langle x\rangle^{2+\delta}\hat f_n(x)dx\leq C\right)\\
		&&\bigskip \cap\left(\int\langle u\rangle^{2\alpha}|\mathcal F \hat f_n(u)|^2 du\leq C\right).
	\end{eqnarray*}
	We will show below that if $E_n$ holds for all $n\in\mathbb{N}$, then conditionally on the original observations,
	\begin{equation}\label{eq: donskerprocess}
		\sqrt n(\widehat{\vartheta_t^*}-\vartheta_t^*)_{t\in\mathbb{R}} \rightsquigarrow \mathbb{G}_F\ in\ \ell_\infty(\mathbb{R}).
	\end{equation}
	This means that $\mu_n^*$ weakly converges to $\mu$ so that 
	$\lim_{n\to\infty}\pi(\mu_n^*,\mu)\mathds{1}_{E_n}=0$, which implies that $\pi(\mu_n^*,\mu)\mathds{1}_{E_n}$ converges in probability to zero. Hence, for all $\epsilon>0$ one has
	\begin{eqnarray*}
		\mathbb{P}(\pi(\mu_n^*,\mu)>\epsilon)&\leq&\mathbb{P}(\pi(\mu_n^*,\mu)\mathds{1}_{E_n}>\epsilon)+(1-\mathbb{P}(E_n))\\
		&=& o(1)+1-\mathbb{P}(E_n).
	\end{eqnarray*}
	Now, it follows from Assumption {\bf B} that for arbitrary $\eta>0$, one can choose $C>0$ large enough so that $\mathbb{P}(E_n)\geq 1-\eta$ for sufficiently large $n$ and therefore,
	\begin{eqnarray*}
		\lim_{n\to\infty}\mathbb{P}(\pi(\mu_n^*,\mu)>\epsilon)&\leq&\eta.
	\end{eqnarray*}
	Letting $\eta\to 0$ proves that for all $\epsilon>0$ one has
	\begin{eqnarray*}
		\lim_{n\to\infty}\mathbb{P}(\pi(\mu_n^*,\mu)>\epsilon)=0,
	\end{eqnarray*}
	which proves the theorem.
	
	It remains to prove \eqref{eq: donskerprocess}, assuming that $E_n$ holds for all $n$ and that $\|\hat f_n-f\Vert_\infty$ converges to zero almost surely. For this task, for an arbitrary fixed $t\in\mathbb{R}$ we decompose the error into a stochastic error term and a bias term as follows
	\begin{eqnarray*}
		\widehat{\vartheta_t^*}-\vartheta_t^*=S_{nt}+B_{nt}
	\end{eqnarray*}
	where
	\begin{eqnarray*}
		S_{nt}=\int\zeta(x-t)\mathcal F^{-1}\left[\mathcal FK_h\frac{\varphi_n^*-\hat \varphi_n}{\varphi_\epsilon}\right](x)dx,
	\end{eqnarray*}
	and
	\begin{eqnarray*}
		B_{nt}=\int\zeta(x-t)(K_h\star \hat f_n-\hat f_n)(x)dx.
	\end{eqnarray*}
	Here, $\hat\varphi_n$ denotes the Fourier transform of the distribution of the $Y_j^*$s. We recall that the distribution of the $\varepsilon_i^*$s is the same as that of the $\varepsilon_i$s whence the bootstrap version of $\varphi_\epsilon$ is $\varphi_\epsilon$ itself and therefore,
	\begin{equation*}
		\hat\varphi_n(u)=\varphi_\epsilon(u)\int e^{iux}\hat f_n(x)dx.
	\end{equation*}
	Since we have assumed that $E_n$ holds for all $n,$ we have the uniform bounds
	\begin{eqnarray*}
		\sup_t|\hat f_n(t)|\leq C, \quad \int\langle u\rangle^{2\alpha}|\mathcal F \hat f_n(u)|^2 du\leq C\quad\mbox{and}\quad \int |x|^2\hat f_n^2(x)dx\leq C,
	\end{eqnarray*}
	so we can conclude as in \cite[Section 4.1.1.]{sohl2012uniform}, where the only difference is that $f_X$ has to be replaced by $\hat f_n$ that is uniformly controlled, that the bias term $B_n$ 
	satisfies
	\begin{eqnarray*}
		\sup_t|B_{nt}|=o(n^{-1/2}).
	\end{eqnarray*}
	Hence, \eqref{eq: donskerprocess} is equivalent to
	\begin{equation}\label{eq: donskerprocessSn}
		(\sqrt nS_{nt})_{t\in\mathbb{R}} \rightsquigarrow\mathbb{G}_F\ in\ \ell_\infty(\mathbb{R}).
	\end{equation}
	To deal with the stochastic term $S_n$ we use the smoothed adjoint inequality (13) in \cite{sohl2012uniform} to obtain that (similar to Equation (15) in that paper) 
	\begin{eqnarray}\label{eq: Snt13}
		S_{nt}=\int\ \mathcal F^{-1}\left[\varphi_\epsilon^{-1}(-\ {\bullet})\mathcal FK_h\right|\star \zeta_t(x)(\mathbb{P}_n^*-\widehat {\mathbb{P}}_n)(\mbox{d} x),
	\end{eqnarray}
	where $\zeta_t(x)=\zeta(t-x)$ for all $x\in\mathbb R$,  $\mathbb{P}_n^*$ is the empirical distribution of the bootstrap sample $Y_1^*,\dots, Y_n^*$, and $\widehat{\mathbb P}_n$ is the common distribution of the bootstrap observations $Y_i^*$s conditionally on the original observations $Y_j$s, that is the distribution with density $\hat f_n\star f_\varepsilon$ with $f_\varepsilon$ being the density of the $\varepsilon_i$s. 
	
	Proving \eqref{eq: donskerprocessSn} taking advantage of the above empirical form for $S_n$  generalizes arguments from the proof of  \cite[Theorem 1]{sohl2012uniform} whence details are omitted and we only point out the main steps and differences (which comes from the fact that in our case, the underlying distribution of the observations depends on $n$). Denoting by $\nu_n^*$ the conditionnal distribution of $\sqrt n S_n$, we will apply  \cite[Theorem 1.5.4.]{aadbook} according to which it suffices to prove that the finite dimensional marginals of $\nu_n^*$ converge to that of $\mu$, and that the sequence $(\nu_n^*)_{n\in\mathbb{N}}$ is asymptotically tight,  in the sense that for all $\epsilon>0$ there exists a compact $K\subset\ell_\infty(\mathbb{R})$ such that
	\begin{equation*}
		\liminf_{n\to\infty}\nu_n^*(K^\eta)\geq 1-\epsilon\quad\mbox{ for all }\eta>0,
	\end{equation*}
	where $K^\eta$ is the set of all functions $f\in\ell_\infty(\mathbb{R})$ such that $\|f-K\|_\infty<\eta$.

	Convergence of  the finite dimensional marginals is given in the Lemma \ref{Lem: fidis} below.  
	\begin{lemma}\label{Lem: fidis}
		Under the assumptions of Theorem \ref{theo: donsker}, if $E_n$ holds for all $n$ and some small enough $\delta>0$, and  $\|\hat f_n-f\Vert_\infty$ converges to zero almost surely as $n\to\infty$, then 
		the finite dimensional distributions of $\sqrt n S_n$ converge towards that of $\mathbb{G}_F$, conditionally on the original observations.
	\end{lemma}

	To prove asymptotic tightness, we split the process into three parts
	\begin{eqnarray*}
		\sqrt n S_n=
		\sqrt n\int(T_1(x)+T_2(x)+T_3(x))({\mathbb{P}}_n^*-\widehat {\mathbb{P}}_n)(\mbox{d} x) 
	\end{eqnarray*}
	where we define $T_1$, $T_2$ and $T_3$ in the same way as in Section 4.2.2 of  \cite{sohl2012uniform}, see their Equations (28), (29) and (30). As mentionned in that paper, due to \cite[Theorem 1.5.7.]{aadbook},  it suffices to show that the three processes $\sqrt n\int T_j 
	( \mbox{d}{\mathbb{P}}_n^*- \mbox{d}\widehat {\mathbb{P}}_n)$ are asymptotically tight, for $j\in\{1,2,3\}$.
	
	In Section 4.2.2 of  \cite{sohl2012uniform}, it is proved that $T_3$, and $T_1$ if $\gamma_s>\beta+1/2$, belong to a bounded subset of $H^{\alpha}(\mathbb{R})$ where $\alpha=1/2+\eta$ for some $\eta>0$ small enough.  In that paper, they conclude from Proposition 1 in \cite{nickl2007bracketing}  that $B$ is  $\mathbb{P}$-Donsker where $\mathbb{P}$ denotes the distribution of the observations they consider. In our case, the underlying distribution $\widehat{\mathbb{P}}_n$ depends on $n$ so we cannot use a similar argument. Note however that proving a Donsker property is much stronger than proving tightness, while only tightness is required here. 
	Nevertheless, asymptotic tightness is proven in Lemma \ref{lem: tightT13} below by extenging part of the arguments from  \cite{marcus1985relationships} to the case of an underlying distribution that depends on $n$. Note that the lemma does not require that the event $E_n$ holds.
	
	\begin{lemma}\label{lem: tightT13} 
		Under the assumptions of Theorem \ref{theo: donsker},  the processes $$\sqrt n\int T_j 
		( \mbox{d}{\mathbb{P}}_n^*- \mbox{d}\widehat {\mathbb{P}}_n)$$ with $j=3$, and $j=1$ if $\gamma_s>\beta+1/2$, are asymptotically tight  conditionally on the original observations. 
	\end{lemma}
	
	To deal with the case $j=2$ we assume that the second event that defines $E_n$ in \eqref{eq: En} is satisfied (see \eqref{eq: hypotight} below). Hence, the result still holds if we assume that $E_n$ holds for all $n$. In Section 4.2.2 of  \cite{sohl2012uniform}, it is proved that the functions $T_2(x)/(1+ix)$ belong to a bounded subset  of $H^{\alpha}(\mathbb{R})$ where $\alpha=1/2+\eta$ for some $\eta>0$ small enough. We combine this result with  empirical process theory and an entropy bound, taken from \cite{nickl2007bracketing}, for the set $H^{\alpha}(\mathbb{R})$, to obtain the following lemma.

	\begin{lemma}\label{lem: tightT2} 
		If Assumption {\bf B} holds,  then the process $\sqrt n\int T_2 
		( \mbox{d}{\mathbb{P}}_n^*- \mbox{d}\widehat {\mathbb{P}}_n)$ is asymptotically tight  conditionally on the original observations.
	\end{lemma}
	
	It remains to deal with the case $j=1$ when $\gamma_s\leq\beta+1/2$. Fix $\xi>0$ and defined
	\begin{eqnarray*}
		K_h^{(0)}=K_h\mathds{1}_{[-\xi,\xi]}.
	\end{eqnarray*}
	In Section 4.2.3 of  \cite{sohl2012uniform}, it is proved that if $\gamma_s> \beta$, and if we define $T$ similarly as $T_1$ with $K_h-K_h^{(0)}$ instead of $K_h$, then $T$ is contained in a bounded subset of $H^{\alpha}(\mathbb{R})$ where $\alpha=\gamma_s-\beta+1>1/2$. Hence, the same arguments as in the proof of Lemma \ref{lem: tightT13}  imply that $\sqrt n\int T
	( \mbox{d}{\mathbb{P}}_n^*- \mbox{d}\widehat {\mathbb{P}}_n)$ is asymptotically tight  conditionally on the original observations. 
	
	It remains to deal with $T_1^{(0)}$, that is defined similarly as $T_1$ with the truncated kernel $K_h^{(0)}$ instead of $K_h$. The arguments in Section 4.2.3 of  \cite{sohl2012uniform} easily generalize to the setting where a fix probability $\mathbb{P}$ is replaced by $\widehat{\mathbb{P}}_n$ if we assume that $\hat f_n (x)\leq C$ for all $x\in\mathbb{R}$. Indeed, our limiting process is the same as that in that paper so we know that it is pre-Gaussian, so in the bootstrap setting it remains  to verify the five conditions  of Theorem 12 in \cite{nickl2012donsker} as in \cite{sohl2012uniform}, but with $\mathbb{P}$ replaced by $\widehat{\mathbb{P}}_n$. Condition (a) does not depend on the underlying probability. To verify condition (b) we write
	\begin{eqnarray*}
		Q_\tau'=\{r-q|r,q\in Q, \left(\int|r-q|^2\mbox{d} \widehat {\mathbb{P}}_n\right)^{1/2}\leq\tau\}
	\end{eqnarray*}
	and we notice that the bounds in  \cite{sohl2012uniform} are uniform with respect to the underlying probability $\mathbb{P}$, whence it holds if the underlying probability is $ \widehat {\mathbb{P}}_n$.
	Similarly, the check of (c)  and (f) relies on bounds that are uniform in the underlying probability. The check of (d) relies on the fact that the density function of the underlying probability is bounded. This condition continues to hold in our setting if we assume that $\hat f_n (x)\leq C$ for all $x\in\mathbb{R}$. Hence, we conclude that  $\sqrt n\int T_1^{(0)}
	( \mbox{d}{\mathbb{P}}_n^*- \mbox{d}\widehat {\mathbb{P}}_n)$ is asymptotically tight  conditionally on the original observations. This completes the proof of Theorem \ref{theo: donsker}. \hfill{$\Box$}

	\subsection{Proof of Lemma \ref{Lem: fidis}}
	We argue conditionally on the original observations.
	It follows from \eqref{eq: Snt13} that for a fixed $t$, $\sqrt n S_{nt}$ is the mean of iid random variables $g_{th}(Y_1^*),\dots,g_{th}(Y_n^*)$ where
	\begin{equation*}
		g_{th}(Y_j^*)=\mathcal F^{-1}\left[(1/\varphi_\epsilon)(-\ {\bullet})\mathcal F  K_h\right]\star \zeta_t(Y_j^*)
	\end{equation*}
	minus expectation. Note that the variance of those variables is
	\begin{eqnarray*}
		v_{nt}
		&=&\int g_{th}^2\mbox{d} \widehat {\mathbb{P}}_n-\left(\int g_{th}\mbox{d} \widehat {\mathbb{P}}_n\right)^2.
	\end{eqnarray*}
	Since $E_n$ holds for all $n$, we have the uniform bounds $\sup_t|\hat f_n(t)|\leq C$, $\|\hat f_n\star f_\epsilon\|_\infty\leq C\int f_\epsilon(x)\mbox{d} x\leq C$ and
	\begin{equation*}
		\int(1+x^2)^{(2+\delta)/2}\widehat{\mathbb{P}}_n(\mbox{d} x)= \int(1+x^2)^{(2+\delta)/2}\hat f_n(x)\mbox{d} x\leq C.
	\end{equation*}
	Using these bounds and similar arguments as for the control of the term in Equation (16) of  \cite{sohl2012uniform}, one can prove that there exists $C'>0$ such that for all $h$ and $n$,
	\begin{equation}\label{eq: lyapounov}
		\int\ \left| g_{th}(x)\right|^{2+\delta} \widehat{\mathbb{P}}_n(\mbox{d} x)\leq C'\Vert\zeta_t\Vert_{Z^{\beta+\delta,1/2+\beta+\delta}}
	\end{equation} 
	(using the same notations as in \cite{sohl2012uniform}), which is finite if $\beta+\delta<1/2$ (see also Example 1 in \cite{sohl2012uniform}). Since $\beta<1/2$, on can choose $\delta$ small enough so that the latter condition holds, hence the bound in \eqref{eq: lyapounov} is finite.

	Now, we consider convergence of $v_{nt}$. It follows from the previous display combined to Jensen's inequality that there exists $C>0$ such that for all integers $n$,
	\begin{equation*}
		\int \left| g_{th}(x)\right|^{2} \hat f_n(x) \mbox{d} x\leq C,
	\end{equation*} 
	which means that the functions $\left| g_{th}\right|^{2} \hat f_n$ all belong to $L_1(\mathbb R)$ with a $L_1$-norm bounded by $C$. This implies that for all $\epsilon>0$ one can find $R>0$ independent of $n$ such that 
	\begin{equation*}
		\sup_n \int_{|x|>R} \left| g_{th}(x)\right|^{2} \hat f_n(x) \mbox{d} x \leq\epsilon,
	\end{equation*}
	see e.g. Theorem 2.32 in \cite{folland1999real}. One can choose $R$ in such a way that, in addition,
	\begin{equation*}
		\int_{|x|>R} \left| g_{t}(x)\right|^{2} f(x) \mbox{d} x \leq\epsilon,
	\end{equation*}
	where $g_t= \mathcal F^{-1}\left[(1/\varphi_\epsilon)(-\ {\bullet})\right]\star \zeta_t$, by integrability of $|g_t|^2f$.
	Now, the functions $\left| g_{th}\right|^{2} \hat f_n$ are bounded on $[-R,R]$: there exists a constant $A>0$ such that the supremum over $n$ of the sup-norms  over $[-R,R]$  of the functions is bounded above by $A$. The constant $A$ is integrable with respect to the Lebesgue measure on $[-R,R]$, so it follows from the dominated convergence theorem under the assumptions that $h$ converges to zero and $\hat f_n$ converges to $f$ that
	\begin{equation*}
		\int_{|x|\leq R} \left| g_{th}(x)\right|^{2} \hat f_n(x) \mbox{d} x-
		\int_{|x|\leq R}\left|g_t(x)\right|^2f(x) \mbox{d} x=o(1).
	\end{equation*}
	Hence, splitting the integrals over $\mathbb R$  into the sum of integrals over $[-R,R]$ and $\mathbb R\backslash [-R,R]$ yields
	\begin{eqnarray*}
		&&\left|\int \left| g_{th}(x)\right|^{2} \hat f_n(x) \mbox{d} x-
		\int\left|g_t(x)\right|^2f(x) \mbox{d} x\right|\\
		&&\qquad\qquad \leq o(1)+ \int_{|x|>R} \left| g_{th}(x)\right|^{2} \hat f_n(x) \mbox{d} x+ \int_{|x|>R} \left| g_{t}(x)\right|^{2} f(x) \mbox{d} x\\
		&&\qquad\qquad \leq o(1)+ 2\epsilon.
	\end{eqnarray*}
	Letting $n\to\infty$ and then $\epsilon\to 0$ one obtains
	\begin{eqnarray*}
		\lim_{n\to\infty}\int \left| g_{th}(x)\right|^{2} \hat f_n(x) \mbox{d} x=
		\int\left|g_t(x)\right|^2f(x) \mbox{d} x.
	\end{eqnarray*}
	Similarly, $\int  g_{th}\hat f_n$ converges to $
	\int g_tf$, 
	so we conclude that $v_{nt}$ converges to 
	$\int\left|g_t\right|^2f-(\int g_tf )^2,$
	which is equal to the variance of the centered Gaussian variable
	$\mathbb G_F(t).$ Combining with \eqref{eq: lyapounov} and the central limit theorem under the Lyapounov condition proves that  $\sqrt n S_{nt}$ converges in distribution to $\mathbb G_F(t).$
	
	The proof can be extended (at the price of cumbersome notations and using a vectoriel central limit theorem) to show that for arbitrary fixed $k\in\mathbb{N}$ and $t_1,\dots,t_k$, the joint distribution of $\sqrt nS_{nt_1},\dots, \sqrt nS_{nt_k}$ converges  to that of $\mathbb G_F(t_1),\dots,\mathbb G_F(t_k)$, which yields  Lemma \ref{Lem: fidis}.
	\hfill{$\Box$}

	\subsection{Proof of Lemma \ref{lem: tightT13}}
	Recall that from Section 4.2.2 of  \cite{sohl2012uniform}, we know that $T_j$ with $j=3$, and $j=1$ if $\gamma_s>\beta+1/2$, belongs to a bounded subset $B$, say, of $H^{\alpha}(\mathbb{R})$ where $\alpha=1/2+\eta$ for some $\eta>0$ small enough. Note that by definition of $T_j$, $B$ can be indexed by $t\in\mathbb{R}$: 
	\begin{equation}\label{eq: Bt}
		B=\{g_t,\ t\in \mathbb{R}\}.
	\end{equation}
	As mentionned  in the paragraph before Theorem 1.3 of that paper, $H^\alpha$ is a Hilbert space for arbitrary $\alpha\geq 0$ if we fix the seminorm so it is a norm, whence we can consider an orthonormal basis $(f_k)_{k\in\mathbb{N}}$. 
	Moreover, it follows from Theorem 1.2 of the same paper that for all $f\in H^\alpha$, $f$ is continuous and bounded  with
	\begin{equation}\label{eq: inclusion}
		\|f\|_\infty\leq \|u^{-1}\|_{L^2}\|f\|_{H^{\alpha}}
	\end{equation}
	where $u(x)=\langle x\rangle=(1+x^2)^{1/2}$. Let $T$ be the inclusion function from $H^\alpha$ to the set of continuous and bounded functions on $\mathbb{R}$, whence $T$ is linear with norm that is less than or equal to $\|k^{-1}\|_{L^2}$. Let $m\in\mathbb{N}$ to be chosen later. One can  obtain Equation (*) on page 325  of \cite{marcus1985relationships} with $P$ replaced by the conditionnal probability $\mathbb P^*$ given $(Y_j)_{j\in\mathbb{N}}$; and then that 
	\begin{eqnarray}\label{eq: Lambdam}
		\sup_n\mathbb{P}^*\left(n^{-1/2}\sup_{f\in B}\left| \sum_{j=1}^n(\delta_{Y_j^*}-\widehat{\mathbb{P}}_n)(f)-\Lambda_m(\delta_{Y_j^*}-\widehat{\mathbb{P}}_n)(f)
		\right| \geq m^{-1}\right)\leq m^{-1}
	\end{eqnarray}
	where for all $g_t:=\sum_ka_{kt}f_k\in B$,
	\begin{equation*}
		n^{-1/2}\sum_{j=1}^n\Lambda_m(\delta_{Y_j^*}-\widehat{\mathbb{P}}_n)(g_t)=
		\sum_{k=1}^{i(m)}a_{kt}Z_{nk}
	\end{equation*}
	for some sufficiently large $i(m)$, with
	\begin{equation*}
		Z_{nk}=n^{-1/2}\sum_{j=1}^n(\delta_{Y_j^*}-\widehat{\mathbb{P}}_n)(f_k).
	\end{equation*}
	Thanks to \eqref{eq: Bt} we can rewrite \eqref{eq: Lambdam} as 
	\begin{eqnarray}\label{eq: Lambdamt}
		\sup_n\mathbb{P}^*\left(\sup_{t\in\mathbb{R}}\left| \mathbb{G}_n^*(t)-\mathbb{G}_n(t)
		\right| \geq m^{-1}\right)\leq m^{-1}
	\end{eqnarray}
	where
	\begin{equation}\label{eq: Gn}
		\mathbb{G}_n^*(t)=n^{-1/2}\sum_{j=1}^n(\delta_{Y_j^*}-\widehat{\mathbb{P}}_n)(g_t),\quad
		\mathbb{G}_n(t)=\sum_{k=1}^{i(m)}a_{kt}Z_{nk}.
	\end{equation}
	For arbitrary $k$ one has
	\begin{eqnarray*}
		\mathbb{E}^*\left(Z_{nk}\right)^2&\leq& 
		\mathbb{E}^*\left(f_k(Y_1^*)\right)^2\\
		&\leq& \|f_k\|_\infty^2\\
		&\leq& \|u^{-1}\|_{L^2}^2
	\end{eqnarray*}
	using \eqref{eq: inclusion} together with the assumption that $f_k\in H^{\alpha}$ has norm one.
	Here $\mathbb{E}^*$ denotes the expectation conditionally on $(Y_j)_{j\in\mathbb{N}}$.
	This implies that for arbitrary $\epsilon>0$ there exists $C_\epsilon>0$ such that
	\begin{equation*}
		\sup_n\mathbb{P}^*\left(|Z_{nk}|>C_\epsilon\right)\leq\epsilon.
	\end{equation*}
	Define $K_k:=\{a_kx,\ |x|\leq C_\epsilon\}$ with $x\in\mathbb{R}$, $a_k$ the function from $\mathbb{R}$ to $\mathbb{R}$ such that $a_k(t)=a_{kt}$. We have  $\sup_t|a_k(t)|\leq \sup_t\| g_t\|_{H^\alpha}<\infty$  since $B$ is bounded whence $a_k\in\ell_\infty(\mathbb{R})$, and $K_k$ is a compact subset of $\ell_\infty(\mathbb{R})$ that satisfies
	\begin{eqnarray*}
		\mathbb{P}^*\left(a_kZ_{nk}\in K_k \right)\geq 1-\epsilon.
	\end{eqnarray*}
	We conclude that the sequence $a_kZ_{nk}$ is asymptotically tight in $\ell_\infty(\mathbb{R})$ for arbitrary $k$ and therefore, 
	it follows from \cite[Theorem 1.5.7.]{aadbook} that $\mathbb{G}_n:=\sum_{k=1}^{i(m)}a_{kt}Z_{nk}$ is asymptotically tight as well in $\ell_\infty(\mathbb{R})$.

	Now, fix $\epsilon>0$ and let $m$ large enough so that $m^{-1}\leq \epsilon/2$. By asymptotic tightness, we can find a compact $K_m$ in the span of $\{a_k,\ k\leq i(m)\}\subset\ell
	_\infty(\mathbb{R})$ such that for all $\delta>0$,
	\begin{equation*}
		\sup_n\mathbb{P}^*\left(\mathbb{G}_n\in K_m^{\delta}\right)\geq 1-\epsilon/2
	\end{equation*}
	and in particular,
	\begin{equation}\label{eq: Ztight}
		\sup_n\mathbb{P}^*\left(\mathbb{G}_n\in K_m^{1/m}\right)\geq 1-\epsilon/2.
	\end{equation}
	With $\tilde K$ a compact subset of $\ell_\infty(\mathbb{R})$ that contains $K_m^{2/m}$ we then have
	\begin{eqnarray*}
		\mathbb{P}^*\left(\mathbb{G}_n^*\not\in \tilde K\right)
		&\leq&\mathbb{P}^*\left(\|\mathbb{G}_n^*- K_m\|_\infty>2/m\right)\\
		&\leq&\mathbb{P}^*\left(\|\mathbb{G}_n- K_m\|_\infty>1/m\right)+ \mathbb{P}^*\left(\|\mathbb{G}_n^*-\mathbb{G}_n\|_\infty>1/m\right)\\
		&\leq&\epsilon/2+m^{-1}\\
		&\leq&\epsilon.
	\end{eqnarray*}
	Hence, $\mathbb{G}_n^*$ is asymptotically tight, where $\mathbb{G}_n^*=\sqrt n\int T_j 
	( \mbox{d}{\mathbb{P}}_n^*- \mbox{d}\widehat {\mathbb{P}}_n)$.
	\hfill{$\Box$}
	
	\subsection{Proof of Lemma \ref{lem: tightT2}}
	
	Thanks to Assumption {\bf B}, we know that for arbitrary $\epsilon>0$, there exists $C_\epsilon>0$  such that for all sufficiently large $n$,   one has
	\begin{equation}\label{eq: hypotight}
		\int |x|^{2+\delta}\hat f_n(x)\mbox{d}x\leq C_\epsilon
	\end{equation}
	with probability larger than $1-\epsilon$. Hence, it suffices to show that for arbitrary $\epsilon>0$, the process $\sqrt n\int T_2 
	( \mbox{d}{\mathbb{P}}_n^*- \mbox{d}\widehat {\mathbb{P}}_n)$ is asymptotically tight  conditionally on the original observations provided that \eqref{eq: hypotight} holds for sufficiently large $n$. For this reason, we assume without loss of generality in the sequel that \eqref{eq: hypotight} holds  where $\epsilon>0$ is arbitrary. 
	
	Recall that from Section 4.2.2 of  \cite{sohl2012uniform}, we know that the functions $T_2(x)/(1+ix)$ belong to a bounded subset of $H^{\alpha}(\mathbb{R})$ where $\alpha=1/2+\eta$ for some $\eta>0$ small enough. Hence, the functions $T_2$ belong to a bounded subset $B$ say, of the set of functions $f$ such that $\|f u^{-1}\|_{H^\alpha}<\infty$ where $u(x)=\langle x\rangle=(1+x^2)^{1/2}$. In the sequel, we assume $\eta$ small enough so that $\eta< \delta/2$ where $\delta$ is taken from Assumption {\bf B}, and we consider some $\gamma\in(\eta,\delta/2)$.
	It follows from \eqref{eq: hypotight} that
	\begin{eqnarray*}
		\int \langle x\rangle^{2+\delta}\mbox{d}\widehat{\mathbb{P}}_n(x)&\lesssim&
		1+\int|x|^{2+\delta}\mbox{d}\widehat{\mathbb{P}}_n(x)\\ &=&
		1+\int|x|^{2+\delta}\hat f_n(x)\mbox{d}x\\
		&\leq&1+C_\epsilon.
	\end{eqnarray*}
	Therefore, it follows from H\"older's inequality that
	\begin{eqnarray*}
		\sup_n \int \langle x\rangle^{2+2\gamma}\mbox{d}\widehat{\mathbb{P}}_n(x)\leq (1+C_\epsilon)^{(2+2\gamma)/(2+\delta)}.
	\end{eqnarray*}
	Since $H^\alpha=B_{2,2}^\alpha$, see e.g. the Appendix in   \cite{sohl2012uniform}, Corollary 1 in \cite{nickl2007bracketing} (with $p=q=2$, $\beta=-1$,  and $s=\alpha$, $r=2$) provides a bound for the bracketing metric entropy of $B$: 
	\begin{eqnarray*}
		\sup_n H_{[]}(\delta, B, \|\ .\ \|_{L_2(\widehat{\mathbb{P}}_n)})\lesssim \delta^{-1/\alpha}. 
	\end{eqnarray*}
	Moreover, the third assertion in Proposition 3 of the same paper implies that there exists a real number $K>0$ such that $\sup_{f\in B}|f|\leq K u$, whence $Ku$ is an envelope of $B$. Due to H\"older's inequality, the enveloppe satisfies
	\begin{eqnarray*}
		\sup_n\int (Ku)^2\widehat{\mathbb{P}}_n(x)&=&\sup_n K^2\int \langle x\rangle^{2}\mbox{d}\widehat{\mathbb{P}}_n(x)<\infty,
	\end{eqnarray*}
	so Theorem 2.14.1 in \cite{aadbook} yields that if \eqref{eq: hypotight} holds for sufficiently large $n$ then
	\begin{eqnarray*}
		\sup_n\mathbb{E}^*\left \|\sqrt n\int T_2
		( \mbox{d}{\mathbb{P}}_n^*- \mbox{d}\widehat {\mathbb{P}}_n)\right\|_\infty
		&\lesssim&
		\sup_n\mathbb{E}^*\sup_{f\in B}\left |\sqrt n\int f
		( \mbox{d}{\mathbb{P}}_n^*- \mbox{d}\widehat {\mathbb{P}}_n)\right|\\
		&\lesssim&\int_0^1\sqrt{1+\epsilon^{-1/\alpha}}\mbox{d}\epsilon\\
		&\lesssim&1+\int_0^1\delta^{-1/(2\alpha)}\mbox{d}\delta
	\end{eqnarray*}
	which is finite since $\alpha>1/2$. The lemma follows.
	\hfill{$\Box$}
	
	\bibliographystyle{plain}
	\bibliography{ConcaveCDF}

\end{document}